\def\Title{Transfer principle in quantum set theory}
\def\Author{Masanao Ozawa}

\documentclass[12pt]{article}
\evensidemargin=\oddsidemargin
  \newcommand{\beq}{\begin{equation}}
  \newcommand{\eeq}{\end{equation}}
  \newcommand{\beql}[1]{\begin{equation}\label{eq:#1}}

  \newcommand{\beqa}{\begin{eqnarray}}
  \newcommand{\eeqa}{\end{eqnarray}}
  \newcommand{\beqas}{\begin{eqnarray*}}
  \newcommand{\eeqas}{\end{eqnarray*}}
  \newtheorem{Theorem}{Theorem}[section]
  \newtheorem{Proposition}[Theorem]{Proposition}
  
  \newtheorem{Corollary}[Theorem]{Corollary}

  \newenvironment{Proof}{\begin{trivlist}
    \item[\hskip \labelsep {\em \indent Proof.}]}{\qed\end{trivlist}}
 \newcommand{\qed}{{\em QED}}

  \newcommand{\Q}{{\bf Q}}
  \newcommand{\R}{{\bf R}}
  \newcommand{\V}{{\bf V}}

  \newcommand{\al}{\alpha}
  \newcommand{\be}{\beta}
  \newcommand{\ch}{\chi}

  \newcommand{\et}{\eta}

  \newcommand{\la}{\lambda}
  \newcommand{\mb}{\mbox}
  
  \newcommand{\om}{\omega}
  \newcommand{\ph}{\phi}
 \newcommand{\ps}{\psi}

  \newcommand{\De}{\Delta}
  \newcommand{\Ga}{\Gamma}

  \newcommand{\And}{\wedge}

  \newcommand{\IFF}{\Leftrightarrow}
  
  \newcommand{\Iff}{\Leftrightarrow}
  
  \newcommand{\Inf}{\bigwedge}

  \newcommand{\Not}{\neg}

  \newcommand{\Or}{\vee}

  \newcommand{\RB}{\R^{(\B)}}

  \newcommand{\Sup}{\bigvee}
  \newcommand{\THEN}{\Rightarrow}
  \newcommand{\Then}{\Rightarrow}

  \newcommand{\VB}{\V^{(\B)}}

\renewcommand{\V}{V}

  \newcommand{\beqan}{\begin{eqnarray*}}
  \newcommand{\beqar}[1]{\begin{equation}\label{#1}\begin{array}{l}}

  \newcommand{\dom}{\mbox{\rm dom}}
  
  \newcommand{\eeqar}{\end{array}\end{equation}}

\renewcommand{\iff}{\quad\mbox{iff}\quad}

\newcommand{\bracket}[1]{\langle#1\rangle}
  \newcommand{\rank}{\mbox{\rm rank}}

\newcommand{\bmat}{\left[\begin{array}{rr}}
\newcommand{\emat}{\end{array}\right]}
\newcommand{\bvec}{\left[\begin{array}{r}}
\newcommand{\evec}{\end{array}\right]}

  \newcommand{\bS}{{\bf S}}

  \newcommand{\cA}{{\cal A}}
  \newcommand{\cB}{{\cal B}}
  \newcommand{\cC}{{\cal C}}
  \newcommand{\cD}{{\cal D}}

  \newcommand{\cH}{{\cal H}}

  \newcommand{\cL}{{\cal L}}
  \newcommand{\cM}{{\cal M}}
  \newcommand{\cN}{{\cal N}}
  
  \newcommand{\cP}{{\cal P}}
  \newcommand{\cQ}{{\cal Q}}
  \newcommand{\cR}{{\cal R}}
  \newcommand{\cS}{{\cal S}}

  \newcommand{\tA}{\tilde{A}}
  \newcommand{\tB}{\tilde{B}}

  \newcommand{\tI}{\tilde{I}}
  \newcommand{\tJ}{\tilde{J}}

\newcommand{\VQH}{\V^{(\cQ(\cH))}}
\newcommand{\VQ}{\V^{(\cQ)}}
\newcommand{\VL}{\V^{(\cQ)}}
\renewcommand{\VB}{\V^{(\cB)}}
\newcommand{\RQ}{\R^{(\cQ)}}
\renewcommand{\RB}{\R^{(\cB)}}

\renewcommand{\dom}{\cD}
\renewcommand{\subset}{\subseteq}
\renewcommand{\inf}{\bigwedge}
\renewcommand{\sup}{\bigvee}
\renewcommand{\Then}{\rightarrow}
\renewcommand{\Iff}{\leftrightarrow}

\newcommand{\iin}{\in\cR}
\newcommand{\val}[1]{[\![#1]\!]}
\newcommand{\vval}[1]{[\![#1]\!]_{\cQ}}
\newcommand{\cdom}{\perp\!\!\!\perp\!\!} 
\newcommand{\commutes}{\ {}^{|}\!\!{}_{\circ}\ }
\newcommand{\cuniv}{\underline{\Or}}

\newcommand{\cx}{\check{x}}
\newcommand{\cy}{\check{y}}
\newcommand{\cz}{\check{z}}

\newcommand{\af}{\,\et\,}

  \title{{\bf \Title}}
  \author{\sc\Author \\
  \it\small Graduate School of Information Sciences\\
\it\small T\^{o}hoku University, Aoba-ku, Sendai,  980-8579,
Japan}
  \date{}
  
\begin{document}
\maketitle
\begin{abstract}
In 1981, Takeuti introduced quantum set theory as the quantum counterpart of 
Boolean valued models of set theory by constructing a model of set theory
based on quantum logic represented by the lattice of closed subspaces in a Hilbert
space and showed that appropriate quantum counterparts of ZFC axioms hold
in the model.  Here, Takeuti's formulation is extended to construct a model
of set theory based on the logic represented by the lattice of projections in an arbitrary 
von Neumann algebra.  A transfer principle is established that enables us to transfer
theorems of ZFC to their quantum counterparts holding in the model.  The set of real
numbers in the model is shown to be in one-to-one correspondence with the set of
self-adjoint operators affiliated with the von Neumann algebra generated by
the logic. Despite the difficulty pointed out by Takeuti that equality axioms do not
generally hold in quantum set theory, it is shown that equality axioms hold for any
real numbers in the model.  It is also shown that
any observational proposition in quantum mechanics can be represented 
by a corresponding statement for real numbers in the model with the truth value
consistent with the standard formulation of quantum mechanics, and that the equality
relation between two real numbers in the model is equivalent with the notion
of perfect correlation between corresponding observables (self-adjoint operators)
in quantum mechanics. The paper is concluded with some remarks on the relevance to
quantum set theory of the choice of the implication connective in quantum logic.
\end{abstract}

\section{Introduction}
Since Birkhoff and von Neumann \cite{BN36} introduced quantum logic in 1936
as the semantical structure of propositional calculus for observational propositions 
on a quantum system, 
there have been continued research efforts for introducing
the methods of symbolic logic into the logical structure 
of our recognition on quantum systems.
However, the introduction of basic notions of sets and numbers in quantum logic
was not realized before Takeuti \cite{Ta81} introduced 
quantum set theory in his seminal paper published in 1981.

Quantum set theory has two main origins, quantum logic
introduced by Birkhoff and von Neumann \cite{BN36}
and Boolean valued models of set theory by which
Scott and Solovay \cite{SS67} reformulated
the method of forcing invented by Cohen
\cite{Coh63,Coh66} for the independence proof of the
continuum hypothesis.
Birkhoff and von Neumann \cite{BN36} argued that the
propositional calculus for quantum mechanics is represented
by the lattice of closed linear subspaces of the Hilbert space
of state vectors of the system so that logical operations correspond to lattice
operations of the subspaces such as set intersection (conjunction), 
closure of space sum (disjunction), and orthogonal complement (negation),
whereas the propositional calculus for classical mechanics is represented
by a Boolean logic corresponding to the Boolean algebra of
Borel subsets  of the phase space of the physical system
modulo the sets of Lebesgue measure zero so that logical
operations correspond to set operations in an obvious way.

The correspondence between classical mechanics and
classical logic can be extended nowadays to its 
ultimate form using Boolean valued analysis introduced
by Scott \cite{Sco69a} and developed by Takeuti 
\cite{Ta78,Ta79b,Ta79c,Ta83a,Ta83b,Ta88}
and many followers 
\cite{Eda83,Jec85,KK94,KK99,Nis91,%
83BT,83BH,84CT,85NC,85TN,86BT,90BB,94FN,95SB,%
Smi84} as follows. 
Let $\VB$ be the Boolean valued universe of set theory
constructed from the complete 
Boolean algebra $\cB$ of Borel subsets of
the phase space of a classical mechanical system modulo
the sets with Lebesgue measure zero.
Then, there is a natural correspondence between 
physical quantities of the system and the real numbers
in the model $\VB$, and by the ZFC transfer principle for
Boolean valued models, any theorem on real numbers 
provable in ZFC gives rise to a valid physical statement
on physical quantities of that system.

One of the interesting programs in constructing
quantum set theory is to extend this correspondence
to that between quantum mechanics and quantum set theory.
Let $\VQ$ be the universe of set theory constructed from the complete
orthomodular lattice $\cQ$ of projections in a von Neumann algebra 
$\cM$ of observables of a physical system \cite{BR79}.
Then, it is expected that there is a natural correspondence between
physical quantities (observables) of the system and the real numbers
in the model $\VQ$, and it is also expected that we have a suitable
transfer principle from theorems in ZFC to a valid physical statement on 
physical quantities of that system.

When the von Neumann algebra $\cM$ is abelian, the corresponding 
physical system is considered  as a classical system, and the projection 
lattice
$\cQ$ is a complete Boolean algebra, so that $\VQ$ is nothing but a
Boolean valued model of set theory.

In Ref.~\cite{Ta81}, Takeuti investigated in the case where
the von Neumann  algebra is a type I factor or equivalently the
algebra of all bounded operators on a Hilbert space $\cH$. 
He showed that axioms of ZFC can be transferred to
appropriate valid statements in $\VQ$ and that real numbers in
$\VQ$ naturally correspond to self-adjoint operators on
$\cH$ or equivalently observables of the quantum system
described by $\cH$.

In this paper, we extend Takeuti's investigation to an arbitrary 
von Neumann algebra $\cM$ on a Hilbert space $\cH$.  
Instead of transferring each axiom of ZFC to an appropriate valid
statement in $\VQ$, this paper establishes a unified transfer principle 
from theorems of ZFC to valid statements in $\VQ$,
which enables us to directly obtain transferred theorems 
in $\VQ$ without proving them from the transferred axioms in $\VQ$.
The only restriction to our transfer principle is that a theorem in ZFC 
to transfer should be expressed by a $\De_0$-formula in the language of
set theory; if $\cM$ is abelian, it is known from a fundamental theorem of 
Boolean valued models of set theory that this restriction can be removed.

After establishing the transfer principle, we investigate real numbers
in $\VQ$.   It is rather easy to extend Takeuti's result to show
that there is a one-to-one 
correspondence between real numbers in $\VQ$ and self-adjoint
operators on $\cH$ affiliated with $\cM$ or naturally observables
of the physical system described by $\cM$.  Another important
result about the real  numbers in $\VQ$ concerns the equality
between real numbers.  As previously shown by Takeuti \cite{Ta81}, 
one of the difficult aspects of quantum set theory is that equality
axioms do not generally hold.   Since the equality is one of the most
fundamental relations, the lack of equality axioms has led many
researchers to question the applicability of quantum set theory to a real
science such as quantum mechanics.  Here, we shall show that
the equality axioms generally hold between any real numbers in $\VQ$.  
We also show that the equality between two real numbers is equivalent
to the notion of quantum perfect correlation \cite{05PCN,06QPC}
between two corresponding 
observables in the quantum system.   Thus, despite the above
difficulty of the equality in $\VQ$,  the equality between real numbers 
in $\VQ$ has a clear physical meaning, and thus quantum set theory
is expected to play a crucial role in the future investigations in
the interpretation of quantum mechanics.  

Section \ref{se:QL} collects results on quantum logic used in the
later sections.
Section \ref{se:UQ} introduces the model $\VQ$ of (quantum) set theory 
based on the projection lattice $\cQ$ of a von Neumann algebra $\cM$,  and
proves some basic properties.  Section \ref{se:ZFC} is devoted to proving
the ZFC Transfer Principle that states that any theorem of ZFC
expressed by a $\De_0$-formula holds in the model $\VQ$ 
up to the well-defined truth value in $\cQ$ representing the degree of
commutativity of elements in $\VQ$ appearing in that formula.
Section \ref{se:RN} introduces real numbers in the model $\VQ$, 
and proves that the
equality axioms hold for the reals in  $\VQ$.  
Section \ref{se:AO} proves that the set $\RQ$ of reals in $\VQ$
is in one-to-one correspondence
with the set of self-adjiont operators affiliated with the von Neumann algebra 
$\cM$ generated by the logic $\cQ$, 
and shows that through this correspondence observationally valid propositions on
the physical system is naturally expressed as the valid statements on the reals in $\VQ$. 
This section also proves that the equality relation between two reals in $\VQ$  
is equivalent with the notion of  the quantum perfect correlation 
between the corresponding observables (self-adjoint operators), 
a notion recently introduced into quantum mechanics \cite{05PCN,06QPC}.
Section \ref{se:CR} concludes the present paper with discussions on 
an open problem on the choice of the implication connective.

\section{Quantum logic}
\label{se:QL}

Let $\cH$ be a Hilbert space.
For any subset $S\subseteq\cH$,
we denote by $S^{\perp}$ the orthogonal complement
of $S$, i.e., $S^{\perp}=\{\ps\in\cH|\ 
\mb{$\bracket{\xi,\ps} =0$ for all $\xi\in S$}\}$.
Then, $S^{\perp\perp}$ is the closed linear span of $S$.
Let $\cC(\cH)$ be the set of all closed linear subspaces in
$\cH$. 
With the set inclusion $M\subset N$ as the partial ordering, 
the set $\cC(\cH)$ is a complete
lattice. The lattice operations on $\cC(\cH)$ are characterized by
$M\And N=M\cap N$,
$M\Or N=(M\cup N)^{\perp\perp}$,
$\bigwedge \cS=\bigcap \cS$,
and $\bigvee \cS=(\bigcup \cS)^{\perp\perp}$
for any $\cS\subseteq\cC(\cH)$.
The operation $M\mapsto M^\perp$ 
is  an orthocomplementation
on the lattice $\cC(\cH)$, with which $\cC(\cH)$ is a
complete orthomodular lattice \cite[p.~65]{Kal83}.

Denote by $\cB(\cH)$ the algebra of bounded linear
operators on $\cH$ and $\cQ(\cH)$ the set of 
projections on $\cH$.
We define the {\em operator ordering} on $\cB(\cH)$ by
$A\le B$ iff $\bracket{\ps,A\ps}\le\bracket{\ps,B\ps}$ for
all $\ps\in\cH$. 
For any $A\in\cB(\cH)$, denote by $\cR (A)\in\cC(\cH)$
the closure of the range of $A$, i.e., 
$\cR(A)=(A\cH)^{\perp\perp}$.
For any $M\in\cC(\cH)$,
denote by $\cP (M)\in\cQ(\cH)$ the projection operator 
of $\cH$
onto $M$.
Then, $\cR\cP (M)=M$ for all $M\in\cC(\cH)$
and $\cP\cR (P)=P$ for all $P\in\cQ(\cH)$,
and we have $P\le Q$ if and only if $\cR (P)\subseteq\cR (Q)$
for all $P,Q\in\cQ(\cH)$,
so that $\cQ(\cH)$ with the operator ordering is also a complete 
orthomodular lattice isomorphic to $\cC(\cH)$.
The lattice operations are characterized by 
$P\And Q=\lim_{n\to\infty}(PQ)^{n}$, 
$P^\perp=1-P$ for all $P,Q\in\cQ(\cH)$.

Let $\cA\subseteq\cB(\cH)$.
We denote by $\cA'$ the {\em commutant of 
$\cA$ in $\cB(\cH)$}.
A self-adjoint subalgebra $\cM$ of $\cB(\cH)$ is called a
{\em von Neumann algebra} on $\cH$ iff $1\in\cM$ and
$\cM''=\cM$.
By the double commutation theorem \cite{Tak79},
a self-adjoint algebra $\cM\subseteq\cB(\cH)$ is a von Neumann
algebra if and only if 
$\cM$ is closed under the weak
operator topology, i.e., for all 
nets $A_{\al}\in\cM$ and $A\in\cB(\cH)$, if 
$\bracket{\ps,A_{\al}\ps}\to\bracket{\ps,A\ps}$
for all $\ps\in\cH$ then $A\in\cM$. 
We denote by $\cP(\cM)$ the set of projections in
a von Neumann algebra $\cM$.

We say that $P$ and $Q$ are {\em commuting}, in 
symbols $P\commutes Q$,
iff $[P,Q]=0$, where $[P,Q]=PQ-QP$.
It is well-known that $P\commutes Q$ if and only if 
$P=(P\And Q)\Or(P\And Q^{\perp})$.
All the relations $P\commutes Q$, $Q\commutes P$,
$P^{\perp}\commutes Q$, $P\commutes Q^{\perp}$,
and $P^{\perp}\commutes Q^{\perp}$ are equivalent.
For any subset $\cA\subseteq\cQ(\cH)$,
we denote by $\cA^{!}$ the {\em commutant} 
of $\cA$ in $\cQ(\cH)$, i.e., $\cA^{!}=
\{P\in\cQ(\cH)\mid P\commutes Q \mbox{ for all }
Q\in\cA\}$ \cite[p.~23]{Kal83}.
Then, $\cA^{!}$ is a complete orthomodular sublattice of 
$\cQ(\cH)$, i.e., $\Inf \cS,\Sup \cS,
P^{\perp}\in\cA^{!}$ for any $\cS\subseteq\cA^{!}$ and
$P\in\cA^{!}$.
A {\em logic} on $\cH$ is a subset $\cQ$ of
$\cQ(\cH)$ satisfying $\cQ=\cQ^{!!}$. 
Thus, any logic on $\cH$ is a complete
orthomodular sublattice of $\cQ(\cH)$.
For any subset $\cA\subseteq\cQ(\cH)$, the smallest 
logic including $\cA$ is the logic
$\cA^{!!}$ called the  {\em logic generated by
$\cA$}.

Then, we have

\begin{Proposition}
A subset $\cQ \subseteq\cQ(\cH)$ is a logic on $\cH$ if
and only if $\cQ=\cP(\cM)$ for some von Neumann algebra
$\cM$ on $\cH$.
\end{Proposition}
\begin{Proof}
Let $\cQ\subseteq\cQ(\cH)$ be such that $\cQ^{!!}=\cQ$.
Then, $\cQ'$ is a von Neumann algebra and
obviously $\cQ^!=\cQ'\cap\cQ(\cH)=\cP(\cQ')$.
Thus, we have $\cQ=\cQ^{!!}=(\cQ^{!})'\cap\cQ(\cH)
=(\cP(\cQ'))'\cap\cQ(\cH)=\cQ''\cap\cQ(\cH)=\cP(Q'')$.
 Suppose $\cQ=\cP(\cM)$ for some von Neumann algebra
$\cM$ on $\cH$.
Then, we have $\cM=\cQ''$.
Thus, we have $\cQ^!=\cQ'\cap\cQ(\cH)=\cP(\cQ')$ and
$\cQ^{!!}=[\cQ'\cap\cQ(\cH)]'\cap\cQ(\cH)
=\cQ''\cap\cQ(\cH)=\cP(\cQ'')=\cP(\cM)$.
Thus, we have $\cQ^{!!}=\cQ$.
\end{Proof}

We define the {\em implication} and the 
{\em logical equivalence} on $\cQ$ by  
$P\Then Q=P^{\perp}\Or(P\And Q)$
 and
$P\Iff Q=(P\Then Q)\And(Q\Then P)$.
The following properties are useful in the later
discussions.  For the proofs, see e.g. Hardegree \cite{Har79}.

\begin{Proposition}\label{th:then}
Let $\cQ$ be a logic on $\cH$
and $P_1,P_2,P,Q\in\cQ$.
The following  hold.

(i) $P\Then Q=\cP\{\ps\in\cH\mid P\ps=(P\And Q)\ps\}$.

(ii) $P\Then Q=\cP\{\ps\in\cH\mid P\ps\in\cR(Q)\}$.

(iii) $P\Iff Q=\cP\{\ps\in\cH\mid P\ps=Q\ps\}$.

(iv) $Q\le P_1\Then P_2 $ if and only if $\cP\cR(P_1Q)\le P_2$.

(v) $\cP\cR(PQ)=P\And (P^\perp\Or Q)$.
\end{Proposition}

A logic on $\cH$ is called {\em Boolean}
iff it is a Boolean algebra.

\begin{Proposition}
The following assertions hold for any $\cA\subseteq\cQ(\cH)$.

(i) $\cA$ is  a Boolean logic if and only if 
$\cA=\cA^{!!}\subseteq\cA^{!}$.

(ii) If $\cA\subseteq\cA^{!}$, the subset
$\cA^{!!}$ is the smallest Boolean logic including
$\cA$.
\end{Proposition}
\begin{Proof}
To prove (i), suppose $\cA=\cA^{!!}\subseteq\cA^{!}$.
By the double commutation theorem, $\cA''$ is
the von Neumann algebra generated by $\cA$.
Since $\cA\subseteq\cA^{!}$, the von Neumann algebra
$\cA''$ is abelian, so that $\cA=\cA^{!!}=\cA''\cap\cQ(\cH)$
is a Boolean logic.
Conversely, suppose that $\cA\subseteq\cQ(\cH)$ is a 
Boolean logic.   Then, $P\commutes Q$ for all
$P,Q\in\cA$, so that we have $\cA=\cA^{!!}\subseteq\cA^{!}$,
and assertion (i) follows.
To prove (ii), suppose $\cA\subseteq\cA^{!}$.
Then, $\cA^{!!}$ is a Boolean sublogic including
$\cA$.  Let $\cB$ be a Boolean sublogic including
$\cA$.  Then, we have
$\cA^{!!}\subseteq\cB^{!!}=\cB$.
Thus, $\cA^{!!}$ is the smallest Boolean sublogic including
$\cA$.
\end{Proof}

The following proposition is useful in later discussions.

\begin{Proposition}\label{th:logic}
Let $\cQ$ be a  logic on $\cH$.
The following hold.

(i)  If $P_{\al}\in\cQ$ and 
$P_{\al}\commutes Q$ for all $\al$, then
$(\Sup_{\al}P_{\al})\commutes Q$, 
$\Inf_{\al}P_{\al}\commutes Q$,
and
$Q \And (\Sup_{\al}P_{\al})=\Sup_{\al}(Q\And
P_{\al}).$

(ii) If $P_1,P_2\commutes Q$, then 
$(P_1\Then P_2)\And Q
=[(P_1\And Q)\Then(P_2\And Q)]\And Q$.
\end{Proposition}
\begin{Proof}
For the proof of (i), see Ref.~\cite{Ta81}.  
To prove (ii), let $P_1,P_2\commutes Q$.
By the set theoretical definition of $\And$ and by Proposition \ref{th:then} (ii), 
it suffices to show that 
$P_1\ps\iin (P_2)$ if and only if 
$(P_1\And Q)\ps\iin (P_2\And Q)$ for every $\ps\iin (Q)$.
Let $\ps\iin (Q)$.  
If $(P_1\And Q)\ps\iin (P_2\And Q)$, we have
$P_{1}\ps=P_1Q\ps= (P_1\And Q)\ps\iin( P_{2}\And Q)\subseteq
\cR(P_{2})$, so that $P_1\ps\iin( P_2)$.
Conversely, if $P_1\ps\iin (P_2)$, we have
$(P_1\And Q)\ps=P_1Q\ps=P_1\ps\iin (P_2)$ and
$(P_1\And Q)\ps=QP_{1}\ps\iin (Q)$,
so that $(P_1\And Q)\ps\iin (P_2\And Q)$.
Thus, we have shown that (ii) holds.
\end{Proof}

Let $\cQ$ be a logic on $\cH$.
Let $\cA\subseteq\cQ$.
A {\em Boolean subdomain} of $\cA$ in $\cQ$ is any
$E\in\cA^!\cap\cQ$ such that 
$P_{1}\And E\commutes P_{2}\And E$ for all 
$P_1,P_2\in\cA$.
Denote by $S_{\cQ}(\cA)$ the set of Boolean subdomains of $\cA$ in $\cQ$, i.e., 
\beqas
S_{\cQ}(\cA)=\{E\in\cA^{!}\cap\cQ\mid 
P_{1}\And E\commutes P_{2}\And E
\mb{ for all }P_{1},P_{2}\in\cA\}.
\eeqas
We shall write $S_{\cQ}(P_1,\cdots,P_n)=S_{\cQ}(\{P_1,\cdots,P_n\})$.

For any $E\in S_{\cQ}(\cA)$, the logic generated by
$\cA\And E$ is a Boolean sublogic of $\cQ$,
since $(\cA\And E)''$ is an abelian von Neumann algebra.
The {\em Boolean domain} of $\cA$ in $\cQ$, 
denoted by  $\cdom_{\cQ}(\cA)$,
is defined as the supremum of $S_{\cQ}(\cA)$, i.e., 
\beqas\label{eq:AD}
\cdom_{\cQ}(\cA)=\Sup S_{\cQ}(\cA).
\eeqas
By Proposition \ref{th:logic} (i) we have 
$\cdom_{\cQ}(\cA)\in\cA^{!}\cap\cQ$
and $P_{1}\And \cdom_{\cQ}(\cA)
\commutes P_{2}\And
\cdom_{\cQ}(\cA)$ for all $P_{1},P_{2}\in\cA$, and hence
\beqa
\cdom_{\cQ}(\cA)=\max\{E\in\cA^{!}\cap\cQ\mid 
P_{1}\And E\commutes P_{2}\And E
\mb{ for all }P_{1},P_{2}\in\cA\}.
\eeqa
 See Ref.~\cite[p.~308]{Ta81}
 for the case $\cQ=\cQ(\cH)$. 

\begin{Theorem}\label{th:cuniv}
For any subset $\cA\subseteq\cQ$, we have
$$
\cdom_{\cQ}(\cA)=\cP\{\ps\in\cH\mid
[P_{1},P_{2}]P_{3}\ps=0
\mb{ \rm  for all }P_{1},P_{2},P_{3}\in\cA\},
$$
and $\cdom_{\cQ}(\cA)=\cdom_{\cQ(\cH)}(\cA)$.
\end{Theorem}
\begin{Proof}
Let $F=\max\{E\in\cA^{!}\mid 
P_{1}\And E\commutes P_{2}\And E
\mb{ for all }P_{1},P_{2}\in\cA\}$
and 
$G=\{\ps\in\cH\mid
[P_{1},P_{2}]P_{3}\ps=0
\mb{ \rm  for all }P_{1},P_{2},P_{3}\in\cA\}$.
Let $\ps\iin (F)$ and 
$P_{1},P_{2},P_{3}\in\cA$.
Then,  the operators $P_{1}F$, $P_{2}F$, 
and $P_{3}F$ are mutually commuting.
We have $\ps=F\ps$, and 
$P_{1}P_{2}P_{3}F=P_{1}FP_{2}FP_{3}F=
P_{2}FP_{1}FP_{3}F=P_{2}P_{1}P_{3}F$.
Hence, we obtain
$P_{1}P_{2}P_{3}\ps=P_{2}P_{1}P_{3}\ps$.
Thus, $\ps\in G$ and $\cR(F)\subseteq G$.
Conversely, suppose $\psi\in G$.
Let $M=(\cA\ps)^{\perp\perp}$, the closed linear subspace
spanned by $\cA\ps$.
Then, $M$ is invariant under all
$P\in\cA$, so that $\cP(M)\in\cA^{!}$.
Let $N=\{\xi\in\cH\mid [P_{1},P_{2}]\xi=0\mbox{ for
all }P_1,P_2\in\cA\}$.  Then, $N$
is a closed subspace such that $M\subseteq N$ by
assumption,
so that $[P_{1},P_{2}]\cP(M)=0$.
Since $\cP(M)\in\cA^{!}$, we have $P_{1}\And \cP(M)
\commutes P_{2}\And \cP(M)$.  Thus, we have
$\ps\in M\subseteq \cR(F)$ and $G\subseteq \cR(F)$, 
so that we have shown the relation $\cR(F)=G$.
By the well-known relation 
$\cP\cN(X)=\cP\cR(X^*)^\perp\in\cQ$,
where $\cN(X)=X^{-1}(0)$, for all $X\in\cQ''$,
we have
$\cP( G)=\Inf_{P_1,P_2,P_3}\cP\cN([P_1,P_2]P_3)\in\cQ$.
It follows that $F\in\cQ$, and hence we have
$F=\cdom_{\cQ}(\cA)$.  Therefore, we have proved 
the relation $\cdom_{\cQ}(\cA)=\cP(G)$.  The relation 
$\cdom_{\cQ}(\cA)=\cdom_{\cQ(\cH)}(\cA)$ now follows 
immediately.
\end{Proof}

Henceforth, for any logic $\cQ$ on $\cH$ and any subset $\cA\subseteq\cQ$,
we abbreviate $\cdom(\cA)=$ $\cdom_{\cQ}(\cA)$, which is the common element 
of $\cQ(\cH)$ for any $\cQ$ by the above theorem.
A simpler characterization of the Boolean domain 
is obtained in terms of
von Neumann algebra generated by $\cA$ as follows.

\begin{Theorem}
For any subset $\cA\subseteq\cQ$, we have
$$
\cdom(\cA)=\cP\{\ps\in\cH\mid
[A,B]\ps=0
\mb{ \rm  for all }A,B\in\cA''\}.
$$
\end{Theorem}
\begin{Proof}
Let $E=\cP\{\ps\in\cH\mid
[P_{1},P_{2}]P_{3}\ps=0
\mb{ \rm  for all }P_{1},P_{2},P_{3}\in\cA\}$
and $F=\cP\{\ps\in\cH\mid
[A,B]\ps=0
\mb{ \rm  for all }A,B\in\cA''\}$.
Suppose $\ps\iin (F)$.
Let $P_{1},P_{2},P_{3}\in\cA$.
We have $P_{2}P_{3}\in\cA''$, and hence
$P_{1}(P_{2}P_{3})\ps=
P_{2}P_{3}P_{1}\ps
=P_{2}P_{1}P_{3}\ps$.
It follows that $\ps\iin (E)$.
Since $E=\cdom(\cA)$, we have
$E\in\cA'$.  Let $P,Q\in\cA$.
We have $[P,QE]=[PE,QE]=[P\And E,Q\And E]=0$.
Since $P\in\cA$ was arbitrary, we have $QE\in\cA'$.
Since $Q\in\cA$ was arbitrary, we have $\cA''E\subseteq \cA'$.
Since $\cA'=(\cA'')'$, we
have
$XYE=YXE$ for any
$X,Y\in\cA''$.
Thus, if  $\ps\iin (E)$, then we have $XY\ps=YX\ps$
and $\ps\in\cR(F)$.  Therefore, we conclude $F=\cdom(\cA)$.
\end{Proof}

\section{Universe of quantum sets}
\label{se:UQ}

We denote by $\V$ the universe of sets
which satisfies the Zermelo-Fraenkel set theory
with the axiom of choice (ZFC).
Throughout this paper, 
we fix the language  $\cL(\in)$
for first-order theory with equality 
having a binary relation symbol
$\in$, bounded quantifier symbols $\forall x\in y$,
$\exists x \in y$, and no constant symbols.
For any class $U$, 
the language $\cL(\in,U)$ is the one
obtained by adding a name for each element of $U$.
For convenience, 
we use the same symbol for an element of $U$ and
its name in $\cL(\in,U)$ as well as for the membership
relation and the symbol $\in$.

To each statement $\ph$ of $\cL(\in,U)$, the satisfaction
relation
$\bracket{U,\in} \models \ph$ is defined by the following recursive rules:
\begin{enumerate}\itemsep=0in
\item $ \bracket{U,\in} \models u\in v
\iff u\in v.$
\item $ \bracket{U,\in} \models u = v
\iff u = v.$ 
\item $ \bracket{U,\in} \models \Not \ph \iff  \bracket{U,\in}
\models
\ph 
\mbox{ does not hold}$. 
\item $ \bracket{U,\in} \models \ph_1 \And \ph_2 
\iff \bracket{U,\in}
\models \ph_1 
\mbox{ and } \bracket{U,\in} \models \ph_2$.
\item $ \bracket{U,\in} \models  (\forall x)\,\ph(x) \iff
\bracket{U,\in} \models \ph(u) \mbox{ for all } u \in U$ .
\end{enumerate}
We regard the other logical connectives and quantifiers as
defined symbols.
Our assumption that $\V$ satisfies ZFC
means that if $\ph(x_1,\ldots,x_n)$ is provable in 
ZFC,  i.e., 
$\mb{ZFC}\vdash \ph(x_1,\ldots,x_n)$, then 
$\bracket{\V,\in}\models \ph(u_1,\ldots,u_n)$ for 
any formula $\ph(x_1,\ldots,x_n)$ of $\cL(\in)$ and
all $u_1,\ldots,u_n\in \V$. 

Let $\cQ$ be a logic on $\cH$.
For each ordinal $ {\al}$, let
$$
\V_{\al}^{(\cQ)} = \{u|\ u:\dom(u)\to \cQ\ \mbox{and}
\ \dom(u) \subseteq \bigcup_{\be<\al}V_{\be}^{(\cQ)}\}.
$$
The {\em $\cQ$-valued universe} $\VL$ is defined
by 
$$
  \VL= \bigcup _{{\al}{\in}\mbox{On}} V_{{\al}}^{(\cQ)},
$$
where $\mbox{On}$ is the class of all ordinals. 
It is easy to see that $\cQ_{1}\subseteq\cQ_{2}$ 
if and only if
$\V^{(\cQ_{1})}_{\al}\subseteq\V^{(\cQ_{2})}_{\al}$
for all $\al$.
Thus, every $\VQ$ is a subclass of $\V^{(\cQ(\cH))}$.
For every $u\in\VQ$, the rank of $u$, denoted by
$\rank(u)$,  is defined as the least $\al$ such that
$u\in \VQ_{\al}$.
It is easy to see that if $u\in\dom(v)$ then 
$\rank(u)<\rank(v)$

For $u\in\VQ$, we define the {\em support} 
of $u$, denoted by $L(u)$, by transfinite recursion on the 
rank of $u$ by the relation
$$
L(u)=\bigcup_{x\in\dom(u)}L(x)\cup\{u(x)\mid x\in\dom(u)\}.
$$
For $\cA\subseteq\VQ$ we write 
$L(\cA)=\bigcup_{u\in\cA}L(u)$ and
for $u_1,\ldots,u_n\in\VQ$ we write 
$L(u_1,\ldots,u_n)=L(\{u_1,\ldots,u_n\})$.
Then, we obtain the following characterization of
subuniverses of $V^{(\cQ(\cH))}$.

\begin{Proposition}\label{th:sublogic}
Let $\cQ$ be a logic on $\cH$ and $\al$ an
ordinal. For any $u\in V^{(\cQ(\cH))}$, we have
$u\in\VL_{\al}$  if and only if
$u\in V^{(\cQ(\cH))}_{\al}$ and $L(u)\subseteq\cQ$.  
In particular, $u\in\VL$ if and only if
$u\in\VQH$ and $L(u)\subseteq\cQ$. 
Moreover, $\rank(u)$ is the least $\al$ such 
that $u\in \VQH_{\al}$ for  any $u\in\VL$.
\end{Proposition}
\begin{Proof}Immediate from transfinite induction on
$\al$.
\end{Proof}

Let $\cA\subseteq\VQ$.  The {\em Boolean
domain of $\cA$}, denoted by $\cuniv(\cA)$, is defined by 
$$
\cuniv(\cA)=\cdom L(\cA).
$$
For any $u_1,\ldots,u_n\in\VQ$, we write
$\cuniv(u_1,\ldots,u_n)=\cuniv(\{u_1,\ldots,u_n\})$.

     To each statement $\ph$ of $\cL(\in,\VL)$ 
we assign the
$\cQ$-valued truth value $ \val{\ph}_{\cQ}$ by the following
recursive rules:
\begin{enumerate}\itemsep=0in
\item $\vval{u = v}
= \inf_{u' \in  \cD(u)}(u(u') \Then
\vval{u'  \in v})
\And \inf_{v' \in   \cD(v)}(v(v') 
\Then \vval{v'  \in u})$.
\item $ \vval{u \in v} 
= \sup_{v' \in \cD(v)} (v(v')\And \vval{u =v'})$.
\item $ \vval{\Not\ph} = \vval{\ph}^{\perp}$.
\item $ \vval{\ph_1\And\ph_2} 
= \vval{\ph_{1}} \And \vval{\ph_{2}}$.
\item $ \vval{\ph_1\Or\ph_2} 
= \vval{\ph_{1}} \Or \vval{\ph_{2}}$.
\item $ \vval{\ph_1\Then\ph_2} 
= \vval{\ph_{1}} \Then \vval{\ph_{2}}$.
\item $ \vval{\ph_1\Iff\ph_2} 
= \vval{\ph_{1}} \Iff \vval{\ph_{2}}$.
\item $ \vval{(\forall x\in u)\, {\ph}(x)} 
= \Inf_{u'\in \dom(u)}
(u(u') \Then \vval{\ph(u')})$.
\item $ \vval{(\exists x\in u)\, {\ph}(x)} 
= \Sup_{u'\in \dom(u)}
(u(u') \And \vval{\ph(u')})$.
\item $ \vval{(\forall x)\, {\ph}(x)} 
= \Inf_{u\in \VL}\vval{\ph(u)}$.
\item $ \vval{(\exists x)\, {\ph}(x)} 
= \Sup_{u\in \VL}\vval{\ph(u)}$.
\end{enumerate}

We say that a statement ${\ph}$ of $ \cL(\in,\VL) $
{\em holds} in $\VL$ iff $ \val{{\ph}}_{\cQ} = 1$.
A formula in $\cL(\in)$ is called a {\em
$\De_{0}$-formula}  iff it has no unbounded quantifiers
$\forall x$ or $\exists x$.

\sloppy
\begin{Theorem}[$\De_{0}$-Absoluteness Principle]
\label{th:Absoluteness}
\sloppy  
For any $\De_{0}$-formula 
${\ph} (x_{1},{\ldots}, x_{n}) $ 
of $\cL(\in)$ and $u_{1},{\ldots}, u_{n}\in \VQ$, 
we have
$$
\val{\ph(u_{1},\ldots,u_{n})}_{\cQ}=
\val{\ph(u_{1},\ldots,u_{n})}_{\cQ(\cH)}.
$$
\end{Theorem}
\begin{Proof}
The assertion is proved by the induction on the complexity
of formulas and the rank of elements of $\VQ$.
First, we shall prove by transfinite induction on $\al$ that 
(i) $\vval{u=v}=\val{u=v}_{\cQ(\cH)}$ for any $u,v\in\VQ_{\al}$, and 
(ii) $\vval{u\in v}=\val{u\in v}_{\cQ(\cH)}$ for any $u\in\VQ_{\al}$ 
and $v\in\VQ_{\al+1}$.
If $\al=0$, relations (i) and (ii) trivially hold.
To prove (i), let $u,v\in\VQ_{\al}$.
By induction hypothesis on (ii), we have 
$\vval{u'\in v}=\val{u'\in v}_{\cQ(\cH)}$ and
$\vval{v'\in u}=\val{v'\in u}_{\cQ(\cH)}$
for all $u'\in\dom(u)$ and $v'\in\dom(v)$.
Thus, we have
\beqas
\vval{u=v}
&=&\Inf_{u'\in\dom(u)}(u(u')\Then\vval{u'\in v})
\And
\Inf_{v'\in\dom(v)}(v(v')\Then\vval{v'\in u})\\
&=&\Inf_{u'\in\dom(u)}(u(u')\Then\val{u'\in v}_{\cQ(\cH)})
\And
\Inf_{v'\in\dom(v)}(v(v')\Then\val{v'\in u}_{\cQ(\cH)})\\
&=&
\val{u=v}_{\cQ(\cH)}.
\eeqas
To prove (ii), suppose $u\in\VQ_{\al}$ and $v\in\VQ_{\al+1}$.
If $v'\in\dom(v)$, we have $v'\in\VQ_{\al}$, and hence we have
$\vval{u=v'}=\val{u=v'}_{\cQ(\cH)}$ from the above.
Thus, we have
\beqas
\vval{u\in v}
&=&
\Sup_{v'\in\dom(v)}(v(v')\And\vval{u=v'})\\
&=&
\Sup_{v'\in\dom(v)}(v(v')\And\val{u=v'}_{\cQ(\cH)})\\
&=&
\val{u\in v}_{\cQ(\cH)}.
\eeqas
Therefore, the assertion holds for atomic formulas.
Any induction step adding a logical symbol works
easily, even when bounded quantifiers are concerned,
since the ranges of the supremum and the infimum 
are common for evaluating $\vval{\cdots}$ and 
$\val{\cdots}_{\cQ(\cH)}$. (This would not happen if
we were to consider unbounded quantifiers.)
\end{Proof}

Henceforth, 
for any $\De_{0}$-formula 
${\ph} (x_{1},{\ldots}, x_{n}) $
and $u_1,\ldots,u_n\in\VQ$,
we abbreviate $\val{\ph(u_{1},\ldots,u_{n})}=
\val{\ph(u_{1},\ldots,u_{n})}_{\cQ}$,
which is the common $\cQ$-valued truth value 
in all $\VL$ such  that $u_{1},\ldots,u_{n}\in\VL$.

The universe $\V$  can be embedded in
$\VQ$ by the following operation 
$\vee:v\mapsto\check{v}$ 
defined by the $\in$-recursion: 
for each $v\in\V$, $\check{v} = \{\check{u}|\ u\in v\} 
\times \{1\}$. 
Then we have the following.
\begin{Theorem}[$\De_0$-Elementary Equivalence Principle]
\label{th:2.3.2}
\sloppy  
For any $\De_{0}$-formula 
${\ph} (x_{1},{\ldots}, x_{n}) $ 
of $\cL(\in)$ and $u_{1},{\ldots}, u_{n}\in V$,
we have
$
\bracket{\V,\in}\models  {\ph}(u_{1},{\ldots},u_{n})
\mbox{ if and only if }
\val{\ph(\check{u}_{1},\ldots,\check{u}_{n})}=1.
$
\end{Theorem}
\begin{Proof}
Let ${\bf 2}$ be the sublogic such that ${\bf 2}=\{0,1\}$.
Then, by induction it is easy to see that 
$
\bracket{\V,\in}\models  {\ph}(u_{1},{\ldots},u_{n})
\mbox{ if and only if }
\val{\ph(\check{u}_{1},\ldots,\check{u}_{n})}_{\bf 2}=1
$
for any ${\ph} (x_{1},{\ldots}, x_{n})$ in $\cL(\in)$, 
and this is
equivalent  to $\val{\ph(\check{u}_{1},\ldots,\check{u}_{n})}=1$
for any $\De_{0}$-formula ${\ph} (x_{1},{\ldots}, x_{n})$ 
by the $\De_0$-Absoluteness Principle.
\end{Proof}
 
Takeuti \cite{Ta81} proved that the following modifications
of the equality axioms hold for the case $\cQ=\cQ(\cH)$.

\begin{Theorem}\label{th:equality81}
Let  $\cQ$ be a logic on $\cH$.
For any $u,u',v,v',w\in \VL$, we have the following.

(i) $\val{u=u}=1$.

(ii) $\val{u=v}=\val{v=u}$.

(iii) $\cuniv(u,v,u')\And
\val{u=u'}\And\val{u\in v}\le \val{u'\in v}$.

(iv) $\cuniv(u,v,u')\And\val{u\in v}\And\val{v=v'}
\le \val{u\in v'}$.

(v) $\cuniv(u,v,w)\And\val{u=v}\And\val{v=w}
\le \val{u=w}$.
\end{Theorem}
\begin{Proof}
Takeuti \cite{Ta81} proved the assertions for the case $\cQ=\cQ(\cH)$, and
the assertions for general $\cQ$ follows from the $\De_0$-Absoluteness
Principle.
\end{Proof}

Takeuti \cite{Ta81} gave examples in which 
the transitivity and the substitution laws do not hold without
modifications, so that the relations
\beqas
\val{u=u'}\And\val{u\in v}&\le& \val{u'\in v},\\
\val{u\in v}\And\val{v=v'}&\le& \val{u\in v'},\\
\val{u=v}\And\val{v=w}&\le& \val{u=w}
\eeqas
do not hold in general.

Takeuti \cite{Ta81} introduced $n$-ary relation
symbols $\cuniv(x_0,\ldots, x_n)$ for any $n=2,3,\ldots$
in the language $\cL(\in)$.  We denote by $\cL(\in,\cuniv)$ 
the language $\cL(\in)$ augmented by relation symbols
$\cuniv(x_0,\ldots, x_n)$.
We extend the $\cQ(\cH)$-valued truth value for all the statements
in $\cL(\in,\cuniv,\VQH)$ by the relation
$$
\val{\cuniv(x_0,\ldots, x_n)}
=\cuniv(u_0,\ldots,u_n)
$$
for any $u_0,\ldots,u_n\in\VQ$.

Takeuti \cite{Ta81} showed that the following
axioms and modifications of axioms of the ZFC
holds in $\VQH$:

{\em Axiom of Infinity.}
$\val{\exists x\in \check{\om}(x\in\check{\omega})
\And  \forall x\in\check{\omega}\exists y\in\check{\om}
(x\in y)}=1.$

{\em Axiom of Pair.}
$\cuniv(u,v)\le
\val{\exists x(\cuniv(u,v,x)\And\forall y
(y\in x\Iff y=u\Or y=v)))}.$

{\em Axiom of Union.}
$\cuniv(u)\le
\val{\exists v(\cuniv(u,v)\And
\forall x(\cuniv(x,u)\Then
(x\in v\Iff \exists y\in u(x\in y))))}.$

{\em Axiom of Replacement.}
$
\val{\forall x\in u \exists y \ph(x,y)}
\le
\val{
\exists v \forall x\in u \exists y\in v \ph(x,y)
}.
$

{\em Axiom of Power Set.}
$
\cuniv(u)\le
\val{
\exists v(\cuniv(u,v)\And\forall t
(\cuniv(u,v,t)\Then(t\in v\Iff
\forall x\in t(x\in u))))
 }.
$

{\em Axiom of Foundation.}
$
\cuniv(u)\And \val{\exists x\in u(x\in u)}\le
\val{\exists x\in u\forall y\in u(\Not y\in u)}.
$

{\em Axiom of Choice.}
$
\cuniv(u)\le
\val{\exists v(\cuniv(u,v)\And
\forall x\in u(\exists y\in x\exists ! z\in u (y\in z)
\Then\exists ! y\in x(y\in v)))}.
$

According to the above, Takeuti \cite{Ta81}  concluded
that a reasonable set theory holds in $\VQH$.
However, it is still difficult to say what theorem
holds in $\VQH$, since we have to construct the proof
for  each theorem using the above ``axioms''.
It may be routine, but the above ``axioms'' do not
ensure that we can tell what theorems hold without
constructing proofs.
In the next section, we shall solve this problem by establishing 
a unified transfer of theorems of ZFC to valid statements on $\VQ$.  

\section{ZFC Transfer Principle in Quantum Set Theory}
\label{se:ZFC}

Let $u\in\VQ$ and $p\in\cQ$.
The {\em restriction} $u|_p$ of $u$ to $p$ is defined by
the following transfinite recursion:
\beqas
\dom(u|_p)&=&\{x|_p\mid x\in\dom(u)\},\\
u|_p(x|_p)&=&u(x)\And p
\eeqas
for any $x\in\dom(u)$.
Note that our definition of restriction is simpler than
the corresponding notion given by Takeuti \cite{Ta81},
and we shall develop the theory of restriction 
along with a different line.

\begin{Proposition}\label{th:L-restriction}
For any $\cA\subseteq \VQ$ and $p\in\cQ$, 
we have 
$$
L(\{u|_p\mid u\in\cA\})=L(\cA)\And p.
$$
\end{Proposition}
\begin{Proof}
By induction, it is  easy to see the relation
$
L(u|_p)=L(u)\And p,
$
so that the assertion follows easily.
\end{Proof}
Let $\cA\subseteq\VQ$.  The {\em logic
generated by $\cA$}, denoted by $\cQ(\cA)$, is  define by 
$$
\cQ(\cA)=L(\cA)^{!!}.
$$
For $u_1,\ldots,u_n\in\VQ$, we write
$\cQ(u_1,\ldots,u_n)=\cQ(\{u_1,\ldots,u_n\})$.

\begin{Proposition}\label{th:range}
For any $\De_0$-formula $\ph(x_1,\ldots,x_n)$ in 
$\cL(\in)$ and $u_1,\cdots,u_n\in\VQH$,
we have $\val{\ph(u_1,\ldots,u_n)}\in\cQ(u_1,\ldots,u_n)$.
\end{Proposition}
\begin{Proof}
Let $\cA=\{u_1,\ldots,u_n\}$.
Since $L(\cA)\subseteq\cQ(\cA)$, it follows from
Proposition \ref{th:sublogic} that $u_1,\ldots,u_n\in
V^{\cQ(\cA)}$.
By the $\De_0$-Absoluteness
Principle, we have 
$\val{\ph(u_1,\ldots,u_n)}=
\val{\ph(u_1,\ldots,u_n)}_{\cQ(\cA)}\in \cQ(\cA)$.
\end{Proof}

\begin{Proposition}\label{th:commutativity}
For any 
$\De_{0}$-formula ${\ph} (x_{1},{\ldots}, x_{n})$ 
of $\cL(\in)$ and $u_{1},{\ldots}, u_{n}\in\VQH$, if 
$p\in L(u_1,\ldots,u_n)^{!}$, then 
$p\commutes \val{\ph(u_1,\ldots,u_n)}$
and $p\commutes \val{\ph(u_1|_p,\ldots,u_n|_p)}$.
\end{Proposition}
\begin{Proof}
Let $u_{1},{\ldots}, u_{n}\in\VQ$.
If $p\in L(u_1,\ldots,u_n)^{!}$, then
$p\in \cQ(u_1,\ldots,u_n)^{!}$.  From Proposition 
\ref{th:range},
$\val{\ph(u_1,\ldots,u_n)}\in\cQ(u_1,\ldots,u_n)$,
so that $p\commutes \val{\ph(u_1,\ldots,u_n)}$.
From Proposition \ref{th:L-restriction},
$L(u_1|_p,\ldots,u_n|_p)=L(u_1,\ldots,u_n)\And p$,
and hence $p\in L(u_1|_p,\ldots,u_n|_p)^{!}$, so that
$p\commutes \val{\ph(u_1|_p,\ldots,u_n|_p)}$.
\end{Proof}

We define the binary relation $x_1\subseteq x_2$ by
``$x_1\subseteq x_2$''=``$\forall x\in x_1(x\in x_2)$.''
Then, by definition for  any $u,v\in\VQ$ we have
$$
\val{u\subseteq v}=
\Inf_{u'\in\dom(u)}
u(u')\Then \val{u'\in v},
$$
and we have $\val{u=v}=\val{u\subseteq v}
\And\val{v\subseteq u}$.

\begin{Proposition}\label{th:restriction-atom}
For any $u,v\in\VQ$ and $p\in L(u,v)^{!}$, we have
the following relations.

(i) $\val{u|_p\in v|_p}=\val{u\in v}\And p$.

(ii) $\val{u|_p\subseteq v|_p}\And p
=\val{u\subseteq v}\And p$.

(iii) $\val{u|_p= v|_p}\And p =\val{u= v}\And p$.
\end{Proposition}

\begin{Proof}
We shall prove by transfinite induction on $\al$ that 
(i) holds for all $u\in\VQ_{\al}$ and $v\in\VQ_{\al+1}$ and that
(ii) and (iii) hold for all $u,v\in\VQ_{\al}$.
If $\al=0$, the relations trivially hold.
To prove (ii), let $u,v\in\VQ_{\al}$ and $p\in L(u,v)^{!}$.
Let $u'\in\dom(u)$.
Since $L(u,v)^{!}\subseteq L(u',v)^{!}$, we have $p\in L(u',v)^{!}$.
Then, we have $\val{u'|_p\in v|_p}=\val{u'\in v}\And p$
by induction hypothesis on (i).
Thus, we have
\beqas
\val{u|_p\subseteq v|_p}
&=&
\Inf_{u'\in\dom(u|_p)}(u|_p(u')\Then\val{u'\in v|_p})\\
&=&
\Inf_{u'\in\dom(u)}(u|_p(u'|_p)\Then\val{u'|_p\in v|_p})\\
&=&
\Inf_{u'\in\dom(u)}
(u(u')\And p)\Then(\val{u'\in v}\And p).
\eeqas
We  have $p\commutes u(u')$ by
assumption on $p$, and $p\commutes\val{u'\in v}$
by Proposition \ref{th:commutativity},
so that
$p\commutes u(u')\Then\val{u'\in v}$ and
$p\commutes (u(u')\And p)\Then(\val{u'\in v}\And p)$.
From Proposition \ref{th:logic} (ii) we have
$$
p\And[
(u(u')\And p)\Then(\val{u'\in v}\And p)]
=
p\And(u(u')\Then\val{u'\in v}).
$$
Thus, from Proposition \ref{th:logic} (i) we have
\beqas
p\And\val{u|_p\subseteq v|_p}
&=&
p\And\Inf_{u'\in\dom(u)}
(u(u')\And p)\Then(\val{u'\in v}\And p)\\
&=&
\Inf_{u'\in\dom(u)}p\And[
(u(u')\And p)\Then(\val{u'\in v}\And p)]\\
&=&
\Inf_{u'\in\dom(u)}p\And(u(u')\Then\val{u'\in v})\\
&=&
p\And\Inf_{u'\in\dom(u)}(u(u')\Then\val{u'\in v})\\
&=&
p\And\val{u\subseteq v}.
\eeqas
Thus, we have proved relation (ii) for all $u,v\in\VQ_{\al}$.
Relation (iii) for all $u,v\in\VQ_{\al}$ follows easily from relation (ii).
To prove (i), suppose $u\in\VQ_{\al}$, $v\in\VQ_{\al+1}$, and $p\in L(u,v)^{!}$.
Let $v'\in\dom(v)$. 
Since $L(u,v)^{!}\subseteq L(u,v')^{!}$, we have $p\in L(u,v')^{!}$.
By relation (iii) for $u,v\in\VQ_{\al}$ shown above, we have  
$\val{u|_p=v'|_p}\And p=\val{u=v'}\And p$.
By Proposition \ref{th:commutativity}, we have 
$p\commutes \val{u=v'}$, so that
$v(v'), \val{u=v'}\in\{p\}^{!}$, and hence
$p\commutes v(v')\And \val{u=v'}$.
Thus,  we  have
\beqas
\val{u|_p\in v|_p}
&=&\Sup_{v'\in\dom(v|_p)}
v|_p(v')\And\val{u|_p=v'}\\
&=&
\Sup_{v'\in\dom(v)}
v|_p(v'|_p)\And\val{u|_p=v'|_p}\\
&=&
\Sup_{v'\in\dom(v)}
v(v')\And p \And \val{u|_p=v'|_p}\\
&=&
\Sup_{v'\in\dom(v)}
(v(v')\And \val{u=v'}\And p)\\
&=&
\left(\Sup_{v'\in\dom(v)}
v(v')\And \val{u=v'}\right)\And p,
\eeqas  
where the last equality follows from Proposition \ref{th:logic} (i).
Thus, by definition of $\val{u=v}$ we obtain
the relation $\val{u|_p\in v|_p}=\val{u=v}\And p$,
and relation (i) for all $u\in\VQ_{\al}$ and $v\in\VQ_{\al+1}$
has been proved.
Therefore, the assertion follows from transfinite induction on $\al$.
\end{Proof}

\begin{Proposition}\label{th:V-restriction}
For any $\De_{0}$-formula ${\ph} (x_{1},{\ldots}, x_{n})$ 
of $\cL(\in)$ and $u_{1},{\ldots}, u_{n}\in\VQ$, if 
$p\in L(u_1,\ldots,u_n)^{!}$, then 
$\val{\ph(u_1,\ldots,u_n)}\And p=
\val{\ph(u_1|_p,\ldots,u_n|_p)}\And p$.
\end{Proposition}
\begin{Proof}
We prove the assertion by induction on 
the complexity of  ${\ph} (x_{1},{\ldots},x_{n})$.
From Proposition \ref{th:restriction-atom}, the assertion
holds for atomic formulas.
Then, the verification of every induction step follows 
from the fact that the function $a\mapsto a\And p$ of all $a\in
\{p\}^{!}$ preserves all the supremum and infimum and
satisfies $(a\Then b)\And p=[(a\And p)\Then (b\And p)] 
\And p$ from Proposition \ref{th:logic} (ii)
and $a^{\perp}\And p=(a\And p)^{\perp}\And p$
for  all $a,b\in\{p\}^{!}$.
\end{Proof}

Now, we can prove the following.

\begin{Theorem}[ZFC Transfer Principle]
For any $\De_{0}$-formula ${\ph} (x_{1},{\ldots}, x_{n})$ 
of $\cL(\in)$ and $u_{1},{\ldots}, u_{n}\in\VQ$, if 
${\ph} (x_{1},{\ldots}, x_{n})$ is provable in ZFC, then
we have
$$ 
\cuniv(u_{1},\ldots,u_{n})\le
\val{\ph({u}_{1},\ldots,{u}_{n})}.
$$
\end{Theorem}
\begin{Proof}
Let $p=\cuniv(u_1,\ldots,u_n)$.
Then, we have $a\And p\commutes b\And p$
for any $a,b\in L(u_1,\ldots,u_n)$, and hence 
there is a Boolean sublogic $\cB$ such that 
$L(u_1,\ldots,u_n)\And p\subseteq \cB$.
From Proposition \ref{th:L-restriction},
we have $L(u_1|_p,\ldots,u_n|_p)\subseteq \cB$.
From Proposition \ref{th:sublogic}, we have
$u_1|_p,\ldots,u_n|_p\in \VB$.
By the ZFC Transfer Principle of the Boolean valued
universe \cite[Theorem 1.33]{Bel85}, we have
$\val{\ph(u_1|_p,\ldots,u_n|_p)}_{\cB}=1$. By the
$\De_0$-Absoluteness Principle, we have
$\val{\ph(u_1|_p,\ldots,u_n|_p)}=1$.
From Proposition \ref{th:V-restriction}, we have
$\val{\ph(u_1,\ldots,u_n)}\And p
=\val{\ph(u_1|_p,\ldots,u_n|_p)}\And p
=p$, and the assertion follows.
\end{Proof}

\section{Real numbers in quantum set theory}
\label{se:RN}

Let $\Q$ be the set of rational numbers in $V$.
We define the set of rational numbers in the model $\VQ$
to be $\check{\Q}$.
We define a real number in the model by a Dedekind cut
of the rational numbers. More precisely, we identify
a real number with the upper segment of a Dedekind cut
assuming that the lower segment has no end point.
Therefore, the formal definition of  the predicate $\R(x)$, 
``$x$ is a real number,'' is expressed by
$$
x\subseteq \check{\Q} \And \exists y\in\check{\Q}(y\in
x)
\And \exists y\in\check{\Q}(y\not\in x)\And
\forall y\in\check{\Q}(y\in x\Iff\forall z\in\check{\Q}
(y<z \Then z\in x)),
$$
where ``$x\subseteq \check{\Q}$''=
``$\forall y\in x(y\in\check{\Q})$.''
We define $\R^{(\cQ)}$ to be the interpretation of 
the set $\R$ of real
numbers in $\VQ$ as follows.
$$
\R^{(\cQ)} = \{u\in\VQ|\ \cD(u)=\cD(\check{\Q})
\ \mb{and }\val{\R(u)}=1\}.
$$

\begin{Theorem}
\label{th:RQ}
For any $u\in\R^{(\cQ)}$, we have the following.

(i) $u(\check{r})=\val{\check{r}\in u}$ for all $r\in\Q$.

(ii)  $\cuniv(u)=1$.
\end{Theorem}
\begin{Proof}
Let $u\in\R^{(\cQ)}$ and $x\in\Q$.
Then, we have
$$
\val{\check{x}\in u}
= \Sup_{y\in\Q}(\val{\check{x}=\check{y}}\And u(\check{y}))\\
= u(\check{x}),
$$
and assertion (i) follows.
We have
\beqas
L(u)
&=&\bigcup_{s\in\dom(u)}L(s)\cup\{u(s)\mid s\in\dom(u)\}\\
&=&\bigcup_{s\in\Q}L(\check{s})
\cup
\{u(\check{s})\mid s\in\Q\}\\
&=&\{0,1,u(\check{s})\mid s\in\Q\},
\eeqas
so that it suffices to show that each $u(\check{s})$
with $s\in\Q$ is mutually commuting. 
By definition, we have
$$
\val{\forall y\in\check{\Q}(y\in u\Iff\forall z\in\check{\Q}
(y<z \Then z\in u))}
=1.
$$
Hence, we have
$$
u(\check{s})=\val{\check{s}\in u}=\Inf_{s<t, t\in\Q}\val{\check{t}\in u}.
$$
Thus, if $s_1<s_2$, then
$u(\check{s}_1)\le u(\check{s}_2)$,
so that
$u(\check{s}_1)\commutes u(\check{s}_2)$.
Thus, each $u(\check{s})$
with $s\in\Q$ is mutually commuting,  and assertion (ii) follows.
\end{Proof}

\begin{Proposition}
If $\cuniv(u)\And\val{\R(u)}=1$, then there is a unique
$v\in\R^{(\cQ)}$ such that $\cuniv(u,v)=1$ and 
$\val{u=v}=1$.
\end{Proposition}
\begin{Proof}
Let $\bar{u}\in\VQ$ be such that
$\dom(\bar{u})=\dom(\check{Q})$ and
$\bar{u}(\check{x})=\val{\check{x}\in u}$ for all
$x\in\Q$.
Since $\cuniv(u)=1$, there is a Boolean sublogic $\cB$
such that $L(u)\subseteq L(u)^{!!}=\cB$.
It is easy to see that $\check{Q}, \bar{u}\in \VB$
so that $\cuniv(u,\bar{u})=1$, 
and $\val{\R(u)}_{\cB}=1$ by the $\De_0$-Absoluteness
Principle.  
By definition, we have $\bar{u}(x)\Then\val{x\in u}=1$
for all $x\in\dom(\bar{u})$, and hence  
$\val{\bar{u}\subseteq u}=1$.
On the other hand, 
we have $\val{u\subseteq\check{\Q}}_{\cB}=1$ from
$\val{\R(u)}_{\cB}=1$.  Thus, 
\beqas
\val{u\subseteq \bar{u}}
&=&
\val{u\subseteq \bar{u}}_{\cB}\\
&=&
\val{\forall x(x\in u\Then x\in\bar{u})}_{\cB}\And
\val{u\subseteq\check{\Q}}_{\cB}\\
&=&
\val{\forall x\in\check{Q}(x\in u\Then x\in\bar{u})}_{\cB}\\
&=&
\Inf_{x\in\dom(\check{Q})}
(\val{x\in u}\Then\val{x\in\bar{u}})\\
&=&1.
\eeqas
To show the uniqueness, let $v_1,v_2\in\R^{(\cQ)}$ be
such that $\val{u=v_1}=\val{u=v_2}=1$
and $\cuniv(u,v_1)=\cuniv(u,v_2)=1$.
Let $x\in\Q$. 
Then, 
$\val{\check{x}\in v_1}=
\val{\check{x}\in u}$ follows from
$\val{u=v_1}=\cuniv(u,v_1,\check{x})=1$,
and we have
$\val{\check{x}\in v_2}=\val{\check{x}\in u}$ similarly.
Thus,  we have 
$\val{\check{x}\in v_1}=
\val{\check{x}\in v_2}$.
Since $v_1(\check{x})=\val{\check{x}\in v_{1}}$
and $v_2(\check{x})=\val{\check{x}\in v_{2}}$,
the relation $\val{v_1=v_2}$ follows easily.
\end{Proof}

\begin{Theorem}
If
$(\R(x_{1})\And\cdots\And\R(x_n))\Then
\ps(x_{1},\ldots,x_{n},x_{n+1},\ldots,x_{n+m})$
is a
$\De_{0}$-formula of $\cL(\in)$
provable in ZFC, then 
for any $u_{1},\ldots,u_{n}\in\R^{(\cQ)}$ and 
$u_{n+1},\ldots,u_{n+m}\in \VQ$
we have
$$
\cuniv(u_1,\ldots,u_{n+m})
\le
\val{\ps(u_{1},\ldots,u_{n+m})}.
$$
\end{Theorem}
\begin{Proof}
By the ZFC transfer principle, we have
$$
\cuniv(u_1,\ldots,u_{n+m})
\le
\val{\R(u_{1})\And\cdots\And\R(u_n)}\Then
\val{\ps(u_{1},\ldots,u_{n+m})},
$$
and the assertion follows from 
$\val{\R(u_{1})\And\cdots\And\R(u_n)}=1$ for
any $u_{1},\ldots,u_{n}\in\R^{(\cQ)}$.
\end{Proof}

In what follows, we write $r\And s=\min\{r,s\}$
and $r\Or s=\max\{r,s\}$ for any $r,s\in\R$.
The Boolean domain $\cuniv(u,v)$ of $u,v\in\RQ$ is characterized as
follows.

\begin{Theorem}\label{th:cuniv-pair}
For any $u,v\in\R^{(Q)}$, we have
$$
\cuniv(u,v)=\cP\{\ps\in\cH\mid
[u(\cx),v(\cy)]\ps
=0\mb{ for all }x,y\in\Q\}.
$$
\end{Theorem}
\begin{Proof}
We have $L(u,v)=\{0,1,u(\cx), v(\cx)\mid x\in\Q\}$.
Let $\ps\in\cR\cuniv(u,v)$.  
From Theorem \ref{th:cuniv}, we have 
$[u(\cx),v(\cy)]u(\cz)\ps=0$ for all $x,y,z\in\Q$. 
Taking the limit $z\to\infty$
we have $[u(\cx),v(\cy)]\ps=0$ for all $x,y\in\Q$.
Conversely, suppose $[u(\cx),v(\cy)]\ps
=0$ for all $x,y\in\Q$.
From Theorem \ref{th:cuniv}, it suffices to show that
$[u(\cx),v(\cy)]u(\cz)\ps=0$
and $[u(\cx),v(\cy)]v(\cz)\ps=0$ for any $x,y,z\in\Q$.
We have
$u(\cx)v(\cy)u(\cz)\ps=u(\cx)u(\cz)v(\cy)\ps=
u(\cx\And\cz)v(\cy)\ps$, and 
$v(\cy)u(\cx)u(\cz)\ps=v(\cy)u(\cx\And\cz)\ps=
u(\cx\And\cz)v(\cy)\ps$, so that
$[u(\cx),v(\cy)]u(\cz)\ps=0$.
Similarly, we also have $[u(\cx),v(\cy)]v(\cz)\ps=0$.
This completes the proof.
\end{Proof}

The $\cQ$-valued equality  $\val{u=v}$ for $u,v\in\RQ$ is characterized
as follows.

\begin{Theorem}
\label{th:equality}
For any $u,v\in\R^{(\cQ)}$ we have
$$
 \val{u=v}=\cP\{\ps\in\cH\mid
\val{\check{r}\in u}\ps=\val{\check{r}\in v}\ps
\mbox{ for all }r\in \Q\}.
$$
\end{Theorem}
\begin{Proof}
From Theorem 5.2 (i) we have
\beqas
\val{u=v}
&=&
\Inf_{r\in \Q}(u(\check{r})\Then\val{\check{r}\in v})\And
\Inf_{r\in \Q}(v(\check{r})\Then\val{\check{r}\in u})\\
&=&
\Inf_{r\in \Q}(\val{\check{r}\in u}\Iff\val{\check{r}\in v}).
\eeqas
From Proposition \ref{th:then} (iii), we have
$$
\val{\check{r}\in u}\Iff\val{\check{r}\in v}
=
\cP\{\ps\in\cH\mid \val{\check{r}\in u}\ps=
\val{\check{r}\in v}\ps\}
$$
Thus, the assertion follows easily.
\end{Proof}

\begin{Theorem}\label{th:forcing_real}
For any $u,v\in\R^{(\cQ)}$ 
and $\ps\in\cH$, the following conditions are all
equivalent.

(i)  $\ps\in\cR\val{u=v}$.

(ii) $u(\check{x})\ps=v(\check{x})\ps$ for any $x\in\Q$.

(iii) $u(\check{x})v(\check{y})\ps=v(\check{x}\And\cy)\ps$ for any
$x,y\in\Q$.

(iv)
$\bracket{u(\check{x})\ps,v(\check{y})\ps}
=\|v(\check{x}\And\check{y})\ps\|^{2}$
for any $x,y\in\Q$.

\end{Theorem}
\begin{Proof}
The equivalence (i) $\IFF$ (ii) follows from Theorem \ref{th:equality}.
Suppose (ii) holds.
Then, we have
$
u(\cx)v(\cy)\ps=u(\cx)u(\cy)\ps=u(\cx\And \cy)\ps=v(\cx\And \cy)\ps.
$
Thus, the implication (ii) $\THEN$ (iii) holds.
Suppose (iii) holds.
We have 
$
\bracket{u(\cx)\ps,v(\cy)\ps}
=\bracket{\ps,u(\cx)v(\cy)\ps}
=\bracket{\ps,v(\cx\And\cy)\ps}
=\|v(\cx\And\cy)\ps\|^{2},
$
and hence the implication (iii)$\THEN$(iv) holds.
Suppose (iv) holds.
Then, we have $\bracket{u(\cx)\ps,v(\cx)\ps}=\|v(\cx)\ps\|^{2}$ and
$\bracket{v(\cx)\ps,u(\cx)\ps}=\|u(\cx)\ps\|^2$.  Consequently, we have 
$
\|u(\cx)\ps-v(\cx)\ps\|^{2}
=
\|u(\cx)\ps\|^2+\|v(\cx)\ps\|^{2}-\bracket{u(\cx)\ps,v(\cx)\ps}
-\bracket{v(\cx)\ps,u(\cx)\ps}
=0, 
$
and hence $u(\cx)\ps=v(\cx)\ps$.  Thus, the implication (iv)$\THEN$(ii)
holds, and the proof is completed.
\end{Proof}

The set $\R_\cQ$ of real numbers in $\VQ$ is defined by
\beqas
\R_\cQ=\R^{(\cQ)}\times\{1\}.
\eeqas
The following theorem shows that the equality is an equivalence
relation between real numbers in $\VQ$.

\begin{Theorem}\label{th:equality-axioms}
The following relations hold in $\VQ$.

(i) $\val{(\forall u\in\R_\cQ) u=u}=1$.

(ii) $\val{(\forall u,v\in\R_\cQ) u=v\Then v=u}=1$.

(iii) $\val{(\forall u,v,w\in\R_\cQ) u=v\And v=w\Then u=w}=1$.

(iv) $\val{(\forall v\in\R_\cQ)(\forall x,y\in v)x=y\And x\in v\Then
y\in v}=1$.

(v) $\val{(\forall u,v\in\R_\cQ)(\forall x\in u)x\in u\And u=v\Then x\in
v}=1.$
\end{Theorem}

\begin{Proof}
Relations (i) and (ii) follow from Theorem \ref{th:equality81}.
To prove (iii), let $u,v,w\in\cD(\R_\cQ)=\R^{(\cQ)}$.
Suppose $\ps\in\cR(\val{u=v}\And\val{v=w})$.
Let $r\in \Q$.  Then, we have $\val{\check{r}\in u}\ps
=\val{\check{r}\in v}\ps$  and 
$\val{\check{r}\in v}\ps=\val{\check{r}\in w}\ps$ from Theorem
\ref{th:equality},
and hence $\val{\check{r}\in u}\ps=\val{\check{r}\in w}\ps$.
Since $r$ was arbitrary, we obtain $\ps\in\cR\val{u=w}$.
Thus, we have
$$
\Inf_{u,v,w\in\cD(\R_\cQ)}\val{u=v\And v=w\Then u=w}=1,
$$
and relation (iii) holds. 
To prove (iv), let $v\in\RQ$, and let $x,y\in\cD(v)$.
Then, we have $s,t\in\Q$ such that $x=\check{s}$
and $y=\check{t}$. If $s\not=t$, then $\val{x=y}=0$, and
the relation trivially holds.
If $s=t$, then we have $\val{x\in v}=
\val{\check{s}\in v}=\val{\check{t}\in v}
=\val{y\in v}$, and hence (iv) holds.
To probe (v), let $u,v\in\RQ$, and let $x\in\cD(u)$.
Suppose $\ps\in\cR(\val{x\in u}\And\val{u=v})$.
Let $r\in\Q$ such that $x=\check{r}$.
Then, $\val{\check{r}\in u}\ps=\ps$ and 
$\val{\check{r}\in u}\ps=\val{\check{r}\in v}\ps$ from Theorem 
\ref{th:equality},
and hence $\val{x\in v}\ps=\ps$, so that
$\ps\in\cR\val{x\in v}$.  Thus,  (v) holds.
\end{Proof}

The following theorem shows 
commutativity follows from equality in $\R^{(\cQ)}$.

\begin{Theorem}\label{th:eq-com}
For any $u_1,\ldots,u_n\in\R^{(\cQ)}$, we have
$$
\val{u_1=u_2\And\cdots\And u_{n-1}=u_n}\le 
\cuniv(u_1,\ldots,u_n).
$$
\end{Theorem}
\begin{Proof}
We have $L(u_1,\ldots,u_n)=\{0,1,u_1(\cx),\ldots,u_n(\cx)|\
x\in\cQ\}$. 
Let $\ps\in\cR\val{u_1=u_2\And\cdots\And u_{n-1}=u_n}$.
From Theorem \ref{th:cuniv}, it suffices to show that
$u_i(\cx)u_j(\cy)u_k(\cz)\ps=u_j(\cy)u_i(\cx)u_k(\cz)\ps$ for any
$i,j,k=1,\ldots,n$ and $x,y,z\in\Q$. 
From Theorem \ref{th:equality-axioms} (iii), we have
$\ps\in\cR\val{u_j=u_k}\cap\cR\val{u_i=u_k}$.
From Theorem \ref{th:forcing_real}, we have
$u_i(\cx)u_j(\cy)u_k(\cz)\ps=u_i(\cx)u_k(\cy\And\cz)\ps
=u_k(\cx\And\cy\And\cz)\ps$ and 
$u_j(\cy)u_i(\cx)u_k(\cz)\ps=u_j(\cy)u_k(\cx\And\cz)\ps
=u_k(\cx\And\cy\And\cz)\ps$.
Thus, the above relation follows easily.
\end{Proof}

The following theorem shows that the equality between real numbers
in $\VQ$ satisfies the substitution law for $\Delta_0$-formulas.

\begin{Theorem}[$\De_0$-Substitution Law]\label{th:sub}
For any $\De_{0}$-formula $\ph(x_{1},\ldots,x_{n})$ in $\cL(\in)$, 
we
have
\beqas
\lefteqn{
[\![
(\forall u_{1},\ldots,u_{n},v_{1},
\ldots,v_{n}\in\R_{\cQ})
}\qquad\\
& & \qquad
(u_{1}=v_{1}\And\cdots\And u_{n}=v_{n})\Then
(\ph(u_{1},\ldots,u_{n})\Iff\ph(v_{1},\ldots,v_{n}))]\!]=1
\eeqas
\end{Theorem}
\begin{Proof}
Let $u_{1},\ldots,u_{n},v_{1},
\ldots,v_{n}\in\R^{(\cQ)}$.
From the ZFC Transfer Principle, we have
$$
\cuniv(u_{1},\ldots,u_{n},v_{1},
\ldots,v_{n})\le
\val{
u_{1}=v_{1}\And\cdots\And u_{n}=v_{n}}\Then
\val{
\ph(u_{1},\ldots,u_{n})\Iff\ph(v_{1},\ldots,v_{n})}.
$$
Since $\cuniv(u_{1},\ldots,u_{n},v_{1},
\ldots,v_{n})\in 
L(u_{1},v_{1},\ldots,u_{n},v_{n})^!$, we have
$$
\val{u_{1}=v_{1}\And\cdots\And u_{n}=v_{n}}
\And\cuniv(u_{1},\ldots,u_{n},v_{1},
\ldots,v_{n})
\le
\val{
\ph(u_{1},\ldots,u_{n})\Iff\ph(v_{1},\ldots,v_{n})}.
$$
From Theorem \ref{th:eq-com}, we have
$$
\val{u_{1}=v_{1}\And\cdots\And u_{n}=v_{n}}
=
\val{u_{1}=v_{1}\And\cdots\And u_{n}=v_{n}}
\And\cuniv(u_{1},\ldots,u_{n},v_{1},
\ldots,v_{n}).
$$
Thus, we have
$$
\Inf_{u_{1},\ldots,u_{n},v_{1},\ldots,v_{n}\in\cD(\R_{\cQ})}
\val{(u_{1}=v_{1}\And\cdots\And u_{n}=v_{n})\Then
(\ph(u_{1},\ldots,u_{n})\Iff\ph(v_{1},\ldots,v_{n}))}=1,
$$
and the assertion follows.
\end{Proof}

\begin{Corollary}
For any $\De_{0}$-formula $\ph(x_{1},\ldots,x_{n})$ in $\cL(\in)$
and any $u_{1},\ldots,u_{n},v_{1},
\ldots,v_{n}\in\R_{\cQ}$, we have
\beqas
\val{u_{1}=v_{1}\And\cdots\And u_{n}=v_{n}}
\And
\val{\ph(u_{1},\ldots,u_{n})}\le \val{\ph(v_{1},\ldots,v_{n})}.
\eeqas
\end{Corollary}
\begin{Proof}
From Theorem \ref{th:sub}, we have
\beqas
\val{u_{1}=v_{1}\And\cdots\And u_{n}=v_{n}}
\le
\val{\ph(u_{1},\ldots,u_{n})}\Then\val{\ph(v_{1},\ldots,v_{n})}.
\eeqas
Thus, Proposition \ref{th:then} (iv) leads to
\beqas
\val{u_{1}=v_{1}\And\cdots\And u_{n}=v_{n}}
\And\val{\ph(u_{1},\ldots,u_{n})}
\le
\val{\ph(v_{1},\ldots,v_{n})}.
\eeqas
\end{Proof}

For any $u,v\in\RQ$, $u\le v$, $u< v$, and $u<v\le w$ are $\De_0$-formula such that
\beqas
\mbox{``$u\le v$''}&=&\mbox{``$v\subseteq u$''},\\
\mbox{``$u < v$''}&=&\mbox{``$(u\le v) \And \Not( u=v)$''},\\
\mbox{``$u<v\le w$''}&=&\mbox{``$(u<v) \And (v\le w)$''}.
\eeqas

Recall that for any $r\in R$ the embedding $\check{r}\in\VQ$ satisfies
\beqas
\cD(\check{r})=\{\check{x}\mid r\le x\in \Q\}
\quad\mbox{and}\quad
\check{r}(\check{x})=1
\eeqas
for all $x\in\Q$ with $ r\le x$.
In order to make the counter part of $r\in\R$ in $\RQ$,  
for any $r\in R$, we define $\tilde{r}\in\RQ$ by
\beqas
\cD(\tilde{r})=\cD(\check{\Q})\quad\mbox{and}\quad
\tilde{r}(\check{t})=\val{\check{r}\le \check{t}}
\eeqas
for all $t\in\Q$.

\begin{Proposition}\label{th:order}
Let $r\in\Q$, $s,t\in\R$, and $u\in\RQ$.
We have the following relations.

(i) $\val{\check{r}\in\tilde{s}}=\val{\check{s}\le \check{r}}.$

(ii) $\val{\tilde{s}\le\tilde{t}}=\val{\check{s}\le\check{t}}.$

(iii) ${\displaystyle \val{u\le\tilde{t}}=\Inf_{t<x\in\Q}u(\check{x})}.$
\end{Proposition}
\begin{Proof}
We have
\beqas
\val{\check{r}\in\tilde{s}}
&=&
\Sup_{x\in\cD(\tilde{s})}\val{\check{r}=x}\And\tilde{s}(x)
=
\Sup_{x\in\Q}\val{\check{r}=\check{x}}\And\tilde{s}(\check{x})
=
\tilde{s}(\check{r})
=
\val{\check{s}\le \check{r}},
\eeqas
so that (i) holds.  We have
\beqas
\val{\tilde{s}\le\tilde{t}}
&=&
\val{(\forall x\in\tilde{t})x\in\tilde{s}}
=
\Inf_{x\in\Q}\val{\check{t}\le \check{x}}\Then\val{\check{s}\le \check{x}}
=
\Inf_{t\le x\in\Q}\val{\check{s}\le \check{x}}
=
\val{\check{s}\le \check{t}},
\eeqas
and hence (ii) holds.  We have
\beqas
\val{u\le\tilde{t}}
=
\val{(\forall x\in\tilde{t})x\in u}
=
\Inf_{x\in\Q}\val{\check{t}\le \check{x}}\Then u(\check{x})
=
\Inf_{t\le x\in\Q}u(\check{x}),
\eeqas
so that (iii) holds.
\end{Proof}

\section{Applications to operator theory and quantum mechanics}
\label{se:AO}

Let $\cM$ be a von Neumann algebra on a Hilbert
space $\cH$
and let $\cQ$ be a logic of projections in $\cM$.
Then, $\cM=\cQ''$ and every von Neumann algebra arises
in this way from a logic $\cQ$ on $\cH$.
A closed operator $A$ (densely defined) on $\cH$ is
said to be {\em affiliated} with $\cM$, in symbols $A\,\et\,\cM$, 
iff $U^{*}AU=A$ for any unitary operator $U\in\cM'$.
Let $A$ be a self-adjoint operator (densely defined) on $\cH$
and let $A=\int_{\R} \la\, dE^A(\la)$
be its spectral decomposition, where $\{E^{A}(\la)\}_{\la\in\R}$ 
is the resolution of the identity belonging to $A$.
It is well-known that $A\af\cQ''$ if and only if $E^A({\la})\in\cQ$ 
for every $\la\in\R$.
Denote by $\overline{\cM}_{SA}$ the set of self-adjoint operators 
affiliated with $\cM$.
Two self-adjoint operators $A$ and $B$ are said to {\em commute},
in symbols $A\commutes B$,
iff $E^A(\la)\commutes E^B(\la')$ for every pair 
$\la,\la'$ of reals.

Let $\cB$ be a Boolean logic on $\cH$. 
Takeuti \cite{Ta78} showed that there is a one-to-one correspondence
between $\RB$ and $\overline{(\cB'')}_{SA}$ as follows.
Let $u\in\RB$.
Then, we have $u(\check{r})\in\cB$ for all $r\in\Q$ and the following
are easily checked.

(i) $\displaystyle\Inf_{r\in\Q}u(\check{r})=0.$

(ii) $\displaystyle\Sup_{r\in\Q}u(\check{r})=1.$

(iii) $\displaystyle u(\check{r})=\Inf_{r<s\in\Q}u(\check{s})$ for every $r\in\Q$.\\
In fact, (i) follows from $\val{\exists y\in\check{\Q}(y\not\in u)}=1$,
(ii) follows from $\val{\exists y\in\check{\Q}(y\in
u)}=1$, and (iii) follows from 
$\val{\forall y\in\check{\Q}(y\in u\Iff\forall z\in\check{\Q}
(y<z \Then z\in u))}=1.$

For any $u\in\R^{(\cB)}$ and $\la\in\R$, we define $E^{u}(\la)$ by 
$$
E^{u}(\la)=\Inf_{\la<r\in\Q}u(\check{r}).
$$
Then, we have the following.

(i) $\displaystyle\Inf_{r\in\Q}E^{u}(\la)=0.$

(ii) $\displaystyle\Sup_{r\in\Q}E^{u}(\la)=1.$

(iii) $\displaystyle E^{u}(\la)=\Inf_{\la<\mu}E^{u}(\mu)$ for every
$\la\in\Q$.\\
The above relations show that $\{E^u(\la)\}_{\la\in\R}$ is a resolution of
the identity in $\cB$ and hence by the spectral theorem there
is a self-adjoint operator $\hat{u}\af\cB''$ uniquely
satisfying $\hat{u}=\int_{\R}\la\, dE^u(\la)$.  On the other
hand, let $A\af\cB''$ be a self-adjoint operator. We define $\tilde{A}\in\VB$ by
$$
\dom(\tilde{A})=\dom(\check{\Q})\mbox{ and }
\tilde{A}(\check{r})=E^{A}(r)\mbox{ for all }r\in\Q.
$$ 
Then, it is easy to see that $\tilde{A}\in\RB$ and we have
$(\hat{u})\tilde{}=u$ for all $u\in\RB$ and $(\tilde{A})\hat{}=A$
for all $A\in\overline{(\cB'')}_{SA}$.
Therefore, the correspondence
between $\RB$ and $\overline{(\cB'')}_{SA}$ is a one-to-one correspondence.
We call the above correspondence the {\em Takeuti correspondence}.
It should be noted that for any $a\in\R$, the real $\tilde{a}\in\RQ$ 
corresponds to the scalar operator $a1$ under the Takeuti correspondence.

Now, we have the following.

\begin{Theorem}
Let $\cQ$ be a logic on $\cH$.  The relations 

(i) ${\displaystyle E^{A}(\la)=\Inf_{\la<r\in\Q}u(\check{r})}$ for all $\la\in\Q$,

(ii) $u(\check{r})=E^{A}(r)$ for all $r\in\Q$,\\
for all $u\in\RQ$ and $A\in \overline{(\cQ'')}_{SA}$ 
sets up a one-to-one correspondence between $\RQ$ and $ \overline{(\cQ'')}_{SA}$.
\end{Theorem}
\begin{Proof}
Let $u\in\RQ$.
From Theorem \ref{th:RQ} we have $\cuniv(u)=1$.
Thus, the logic $L(u)^{!!}$ generated by $L(u)$ 
is a Boolean logic.  Let $\cB=L(u)^{!!}$.  By the 
$\De_0$-Absoluteness Principle, $u\in\RB$.
Thus, by Takeuti's result above there is a self-adjoint
operator $A=\hat{u}\af \cB''$ satisfying relations (i) and (ii).
Since $\cB''\subseteq\cQ''$, we have shown that for any
$u\in\RQ$ there exists $A\af \cQ''$ satisfying (i) and (ii),
and the uniqueness of such $A$ follows easily.
On the other hand, let $A\af\cQ''$.  Let $\cB$ be
the Boolean logic generated by $\{E^{A}(\la)\}_{\la\in\R}$.
Then, by Takeuti's result above we have $u=\tilde{A}\in\VB$
satisfying relations (i) and (ii).  Since $A\af\cQ''$, we have 
$\cB\subseteq\cP(\cQ'')=\cQ$, and hence $u\in\VQ$.  
By the $\De_0$-Absoluteness
Principle,  we also have $u\in\RQ$.  Thus, we have proved that
relations (i) and (ii) determine a one-to-one correspondence
between $\RQ$ and $\overline{(\cQ'')}_{SA}$.
\end{Proof}

Let $E^{A}(\la)$ be the resolution of the 
identity belonging to a self-adjoint operator $A$.
Let $a<b\in\R$.  For the interval $I=(a,b]$, we define 
$$
E^{A}(I)=E^{A}(b)-E^{A}(a),
$$
and we define the corresponding interval 
$\tilde{I}$ of real numbers in $\VQ$ by
\beqas
\cD(\tilde{I})=\RQ
\quad\mbox{and}\quad
\tilde{I}(u)=\val{\tilde{a}<u}\And\val{u\le\tilde{b}}
\eeqas
for all $u\in\RQ$.   

\begin{Theorem}\label{th:statistics}
Let $\cQ$ be a logic on $\cH$.
For any self-adjoint operator $A\et\cQ''$ and any interval $I=(a,b]$, 
we have
$$ 
\val{\tilde{A}\in\tilde{I}}=E^{A}(I). 
$$
\end{Theorem}
\begin{Proof}
Let $u\in\RQ$ and $I\in (a,b]$.  We have
\beqas
\val{u\in\tilde{I}}
&=&
\Sup_{v'\in\RQ}\val{\tilde{a}<v'}\And\val{v'\le\tilde{b}}\And\val{u=v'}.
\eeqas
From the $\De_0$-Substitution Law, we have
$$
\val{\tilde{a}<v'}\And\val{v'\le\tilde{b}}\And\val{u=v'}\le
\val{\tilde{a}<u}\And\val{u\le\tilde{b}}
$$
for any $v'\in\RQ$, so that
$$
\val{u\in\tilde{I}}\le\val{\tilde{a}<u}\And\val{u\le\tilde{b}}.
$$
From
\beqas
\Sup_{v'\in\RQ}\val{\tilde{a}<v'}\And\val{v'\le\tilde{b}}\And\val{u=v'}
&\ge& 
\val{\tilde{a}<u}\And\val{u\le\tilde{b}}\And\val{u=u}\\
&=&
\val{\tilde{a}<u}\And\val{u\le\tilde{b}},
\eeqas
we have
$$
\val{u\in\tilde{I}}\ge\val{\tilde{a}<u}\And\val{u\le\tilde{b}}.
$$
Hence, we have
$$
\val{u\in\tilde{I}}=\val{\tilde{a}<u}\And\val{u\le\tilde{b}}
$$
for any $u\in\RQ$. 
Let $A\et\cQ''$ be a self-adjoint operator.
From Proposition \ref{th:order},  we have
$$
E^{A}(\la)=\Inf_{\la<x\in\Q}\tilde{A}(\check{x})=\val{\tilde{A}\le\tilde{\la}}
$$
for any $\la\in \R$.
Thus, we have
\beqas
E^{A}(I)&=&E^{A}(b)-E^{A}(a)
=
E^{A}(b)\And E^{A}(a)^{\perp}
=
\val{\tilde{A}\le\tilde{b}\And \Not( \tilde{A}\le\tilde{a})}
=
\val{\tilde{a}<\tilde{A}\le \tilde{b}}\\
&=&
\val{\tilde{A}\in\tilde{I}}.
\eeqas
Thus,  the assertion follows.
\end{Proof}

Let $\cQ$ be a logic on $\cH$.  Any unit vector $\ps\in\cH$ is called a (vector)
state of $\cQ$.  We define the probability of any statement $\ph$ in
$\cL(\in,\VQ)$ in a state $\ps$ by
$$
\Pr\{\ph\|\ps\}=\|\val{\ph}\ps\|^{2}.
$$
We say that statement $\ph$ in
$\cL(\in,\VQ)$ {\em holds} in state $\ps$ iff $\Pr\{\ph\|\ps\}=1$, and 
this condition is equivalent to $\ps\in\cR\val{\ph}$.
In what follows, we shall show that 
this probabilistic interpretation of the statements in $\cL(\in,\VQ)$ is consistent
with the standard formulation of quantum mechanics.

In the standard formulation of (non-relativistic) quantum mechanics,
every quantum system $\bS$ corresponds to a Hilbert space $\cH$.
An {\em observable} of $\bS$ is represented by a self-adjoint operator
(densely defined) on $\cH$, and a {\em (vector) state} of $\bS$ is represented by
a unit vector $\ps\in\cH$. 

For any observable $A$, let $E^{A}(\la)$ be the resolution of the 
identity belonging to $A$.
A basic principle of quantum mechanics is formulated as follows \cite[p.~200]{vN55}.
{\em In the state $\psi$ mutually commuting observables $A_1,\ldots, A_n$ 
take values from the respective intervals $I_1,\ldots,I_n$ with the probability 
\beq\label{eq:jp}
\|E^{A_1}(I_1)\cdots E^{A_n}(I_n)\ps\|^2.
\eeq
}

Let $\tilde{A}_1,\ldots,\tilde{A}_n$ be the corresponding elements in $\VQH$,
and $\tilde{I}_1,\ldots,\tilde{I}_n$ the corresponding intervals in $\VQH$.
Then, from Theorem \ref{th:statistics} we have
\beqas\label{eq:jp-qst}
\Pr\{\tA_1\in \tI_1 \And\cdots\And\tA_n\in \tI_n\|\ps\}
&=&
\|\val{\tilde{A}_{1}\in\tilde{I}_1\And\cdots\And\tilde{A}_{n}\in\tilde{I}_n}\ps\|^{2}\\
&=&
\|E^{A_1}(I_1)\cdots E^{A_n}(I_n)\ps\|^2.
\eeqas
Thus, we can restate the basic principle of quantum mechanics as follows.
{\em In the state $\psi$ mutually commuting observables $A_1,\ldots, A_n$ 
take values from the respective intervals $I_1,\ldots,I_n$ with the probability 
\beq\label{eq:jp}
\Pr\{\tA_1\in \tI_1 \And\cdots\And\tA_n\in \tI_n\|\ps\}.
\eeq
}

Therefore, we have shown that there is a natural one-to-one correspondence between
observables of a quantum system described by a Hilbert space $\cH$ and real
numbers in the universe $\VQH$ of quantum sets, so that observational propositions
on the quantum system is naturally expressed as the valid statements on the real numbers
in $\VQH$.

In the conventional interpretation of quantum mechanics \cite{vN55}, 
atomic observational
propositions are restricted to those of the form $\tA\in\tI$ for an observable $A$
and an interval $\tI$ as above.  However, quantum set theory is expected to 
extend the interpretation of quantum mechanics to a more general class of
observational propositions.  Here, we introduce one such extension of the
interpretation.  

For any two commuting observables $A$ and $B$ and any state $\ps$,
we have a joint probability distribution   $\mu^{A,B}_{\ps}$ of $A$ and $B$
in $\psi$, a probability measure on $\R^{2}$ satisfying 
$$
\mu^{A,B}_{\ps}(I\times J)
=
\Pr\{\tA\in \tI\And\tB\in \tJ\|\ps\}
=
\|E^{A}(I)E^{B}(J)\ps\|^{2}
$$
for any intervals $I$ and $J$.
Then, it is natural to consider that $A$ and $B$ have the same
value in state $\ps$ if and only if
$$
\Pr\{\tA\in \tI\And\tB\in \tJ\|\ps\}=0
$$
for any $I,J$ such that $I\cap J=\emptyset$, and moreover this condition is equivalent 
to the following conditions:

(i) $\mu^{A,B}_{\ps}(\{(a,b)\in\R^{2}\mid a=b\})=1.$

(ii) $\mu^{A,B}_{\ps}(\{(a,b)\in\R^{2}\mid a\not=b\})=0.$

(iii) $\mu^{A,B}_{\ps}(I\times J)=\mu^{A,B}_{\ps}((I\cap J)\times\R)=
\mu^{A,B}_{\ps}(\R\times (I\cap J))$ for any intervals $I$ and $J$.\\
Following the classical probability theory, we say in this case 
that observables $A$ and $B$ are perfectly correlated in state $\ps$.
Thus, the notion of perfect correlation is straightforward for
any pair of commuting observables. 
However, the problem of extending this notion to any pair of non-commuting
observables has a non-trivial difficulty, since we have no universal definition of
the joint probability distribution for noncommuting observables.

In the recent investigations \cite{05PCN,06QPC},
we have obtained a satisfactory solution for the above problem.
Here, we shall consider this problem in the light of quantum set theory.
Since $\val{\tA\le \tilde{r}}=\val{\check{r}\in\tA}$ for any 
$r\in\Q$, it is natural to say that $A$ and $B$ have the same value 
in state $\ps$ iff $\ps\in\cR\val{(\forall r\in\check{Q})r\in\tA \Iff r\in\tB}$,
or equivalently iff $\ps\in\cR\val{\tA=\tB}$.

Now, we shall show that the above condition is equivalent to the notion of
perfect correlation formulated in \cite{05PCN,06QPC}.
Let $A$ be an observable.
For any (complex-valued) bounded Borel function $f$ on $\R$, 
we define the observable
$f(A)$ by
$$
f(A)=\int_{\R}f(\la)\,dE^{A}(\la).
$$
We shall denote by $B(\R)$ the space of bounded Borel functions on $\R$.
For any Borel set $\De$ in $\R$, we define $E^{A}(\De)$ by
$E^{A}(\De)=\ch_{\De}(A)$, where $\ch_{\De}$ is a Borel function
on $\R$ defined by $\ch_{\De}(x)=1$ iff $x\in\De$ and $\ch_{\De}(x)=0$ iff
$x\not\in\De$. 
For any pair of observables $A$ and $B$, the {\em joint probability distribution}
of $A$ and $B$ in a state $\ps$ is a probability measure $\mu^{A,B}_{\ps}$
on $\R^{2}$ satisfying 
$$
\mu^{A,B}_{\ps}(\De\times\Ga)=
\bracket{\ps,(E^{A}(\De)\And E^{B}(\Ga))\ps}
$$
for any $\De,\Ga\in\cB(\R)$.
Gudder \cite{Gud68} showed that the joint probability distribution $\mu^{A,B}_{\ps}$
exists if and only if the relation $[E^{A}(\De),E^{B}(\Ga)]\ps=0$ holds for
every $\De,\Ga\in\cB(\R)$.

\begin{Theorem}\label{th:qpc}
For any observables
(self-adjoint operators) $A, B$ on $\cH$ and any state
(unit vector) $\ps\in\cH$, the following conditions are all equivalent.

(i)  $\ps\in\cR\val{\tilde{A}=\tilde{B}}$.

(ii) $E^{A}(r)\ps=E^{B}(r)\ps$ for any $r\in\Q$.

(iii) $f(A)\ps=f(B)\ps$ for all $f\in B(\R)$.

(iv) $\bracket{E^{A}(\De)\ps,E^{B}(\Ga)\ps}=0$
for any $\De,\Ga\in\cB(\R)$ with $\De\cap\Ga=\emptyset$.

(v) There is the joint probability distribution $\mu^{A,B}_{\ps}$ of $A$ and $B$
in $\ps$ satisfying 
$$
\mu^{A,B}_{\ps}(\{(a,b)\in\R^{2}\mid a=b\})=1.
$$

\end{Theorem}
\begin{Proof}
The equivalence (i) $\IFF$ (ii)
follows from Theorem \ref{th:forcing_real}.
Suppose that (ii) holds.  Let $\la\in\R$.  
If $r_1,r_2,\ldots$ be a decreasing sequence of rational numbers
convergent to $\la$, then $E^{A}(r_n)\ps$ and $E^{B}(r_n)\psi$ 
are convergent to $E^{A}(\la)\ps$ and $E^{B}(\la)\ps$, respectively, so that
$E^{A}(\la)\ps=E^{B}(\la)\ps$ for all $\la\in\R$.  Thus, we have 
$$
\bracket{\xi,f(A)\ps}=\int_\R f(\la)\,d\bracket{\xi,E^{A}(\la)\ps}=\int_\R
f(\la)\,d\bracket{\xi,E^{B}(\la)\ps}=\bracket{\xi,f(B)\ps}
$$ 
for all $\xi\in\cH$, and hence we have $f(A)\ps=f(B)\ps$ for all
$f\in B(\R)$.   Thus, the implication (ii) $\THEN$ (iii) holds.
Since condition (ii) is a special case of condition (iii) where $f=\ch_{(-\infty,r]}$,
the implication (iii) $\THEN$ (ii) is trivial, so that the equivalence
(ii) $\IFF$ (iii) follows.  The equivalence of assertions (iii), (iv), and (v) have
been already proved in Ref.~\cite{06QPC}, the proof is completed.
\end{Proof}

Condition (iv) above is adopted as the defining condition for $A$ and $B$ to be
perfectly correlated in $\ps$ because of the simplicity and generality of the formulation.
Condition (v) justifies our nomenclature calling $A$ and $B$ ``perfectly 
correlated.''
By condition (i), quantum logic justifies the assertion that ``perfectly 
correlated'' observables actually have the same value in the given state.
For further properties and applications of the notion of perfect correlation,
we refer the reader to Ref.~\cite{06QPC}.

\section{Concluding Remarks}
\label{se:CR}

In classical logic, the implication connective 
$\Then$ is defined by negation and
disjunction as $P\Then Q=(\Not P)\Or Q$.
In quantum logic several counterparts have been proposed.
Hardegree \cite{Har81} proposed the following requirements for the implication
connective.

(E)    $P\Then Q=1$ if and only if $P\le Q$.

(MP) $P             \And (P\Then Q) \le Q$.

(MT) $Q^{\perp}\And (P\Then Q) \le P^{\perp}$.

(LB) If $P\commutes Q$, then $P\Then Q=P^{\perp}\Or Q$.\\
Then, the work of Kotas \cite{Kot67} can be applied to the problem as to what
complemented-lattice-polynomial definitions of $P\Then Q$ satisfy the above
conditions; there are exactly three possibilities:

(i) $P\Then_1 Q=P^{\perp}\Or (P\And Q).$

(ii) $P\Then_2 Q=(P\Or Q)^{\perp}\Or Q.$

(iii) $P\Then_3 Q=(P\And Q)\Or (P^{\perp}\And Q)\Or (P^{\perp}\And Q^{\perp}).$\\
However, so far we have no general agreement on the choice from the above, 
although the majority view favors definition (i), the so-called Sasaki arrow 
\cite{Urq83}.

In quantum set theory, the truth values of atomic formulas, $\val{u\in v}$ and
$\val{u=v}$, depend crucially on the definition of the implication connective.
Takeuti \cite{Ta81} chose the Sasaki arrow for this and the present work has
followed Takeuti's choice.  

Here, another approach should be also mentioned.
Titani and Kozawa \cite{TK03} developed  ``quantum set theory''  based on the
implication connective  defined by $P\Then Q=1$ if $P\le Q$ and  $P\Then Q=0$
otherwise. Then, their implication connective satisfies (E), (MP), and (MT),  
but does not satisfy (LB).  They successfully showed that there is a one-to-one
correspondence between quantum reals and observables.  
However, their truth values of relations between quantum reals is not
consistent with the standard interpretation of quantum mechanics.
In fact, for any observable $A$ and any real number $a\in\R$,
their truth value of the relation $\tilde{A}\le \tilde{a}$ satisfies 
$\val{\tilde{A}\le \tilde{a}}=1$ if $A\le a1$ and $\val{\tilde{A}\le \tilde{a}}=0$
otherwise.  Moreover, 
their equality between two quantum reals takes only two values so that
$\val{\tA=\tB}=1$ if $A=B$ and $\val{\tA=\tB}=0$ otherwise.
Thus, the Titani-Kozawa approach does not lead to consistency results similar to
Theorem \ref{th:statistics} or Theorem \ref{th:qpc}.

Although the choice of the implication connective affects the truth value 
assignment of any statements in the language of set theory, it is  natural
to expect that not only the implication connective (i) but also the implication 
connectives (ii) and (iii) will lead to consistency results similar
to the above mentioned results.
Thus, it is an interesting open problem to figure out the characteristic
differences in quantum set theories with different choices 
of the implication connective from (i), (ii) and (iii) above.
From the investigations in this line, it is expected to find a reasonable
affirmative answer to the question ``Is a quantum logic a logic?'' \cite{GG71}.

\section*{Acknowledgments}
This work was supported 
by the SCOPE project of MIC of Japan and 
the Grant-in-Aid for Scientific Research of
the JSPS.

\providecommand{\bysame}{\leavevmode\hbox to3em{\hrulefill}\thinspace}

\end{document}